\documentclass[11pt,reqno]{amsart}
\usepackage[textwidth=5cm,textheight=5cm]{geometry}
\newgeometry{textwidth=15cm,textheight=20cm}
\usepackage{mathtools,amsmath,amssymb}
\usepackage{url}
\usepackage[colorlinks=true,citecolor=blue]{hyperref}
\usepackage[numbers,sort&compress]{natbib}
\usepackage{amsfonts}
\usepackage{graphicx}
\usepackage{algorithm}
\usepackage{algorithmic}
\usepackage{xcolor}
\usepackage[utf8]{inputenc}
\usepackage{pifont}

\newcommand{\nothing}{}
 
\newcommand{\cmark}{\ding{51}}%
\newcommand{\xmark}{\ding{55}}%
\newcommand{\la}{\left\langle}
\newcommand{\ra}{\right\rangle} 

\newcommand{\alglabel}{%
  \addtocounter{ALC@line}{-1}
  \refstepcounter{ALC@line}
  \label
}
\newcommand{\one}{\mathbf{1}}
\newcommand{\onen}{\one_n}
 
\newtheorem{theorem}{Theorem}[section]
\newtheorem{lemma}[theorem]{Lemma}

\newtheorem{proposition}[theorem]{Proposition}

\theoremstyle{remark}
\newtheorem{remark}[theorem]{Remark}

\newtheorem{definition}[theorem]{Definition}

\begin{document}

\title[Accelerated Bregman Primal-Dual methods for Optimal Transport]{Accelerated Bregman Primal-Dual methods applied to
  Optimal Transport and Wasserstein Barycenter problems}

\author[A. Chambolle]{Antonin Chambolle}
\thanks{}
\address[A. Chambolle]{CEREMADE, CNRS \& Université Paris-Dauphine, Paris, France} 
\email{ \href{chambolle@ceremade.dauphine.fr}{\nolinkurl{chambolle@ceremade.dauphine.fr}}}

\author[J.P. Contreras]{Juan Pablo Contreras}
\address[J.P. Contreras]{Universidad Adolfo Ibáñez, Facultad de Ingenier\'ia y Ciencias, Santiago \& 
Departamento de Ingeniería Industrial, Universidad Católica del Norte, Antofagasta, Chile }
\email{ \href{mailto:juan.contrerasff@gmail.com}{\nolinkurl{juan.contrerasff@gmail.com}}}
\thanks{}

\begin{abstract}
This paper discusses the efficiency of Hybrid Primal-Dual (HPD) type algorithms to approximate solve discrete Optimal Transport (OT) and Wasserstein Barycenter (WB) problems, with and without entropic regularization. Our first contribution is an analysis showing that these methods yield state-of-the-art convergence rates, both theoretically and practically. Next, we extend the HPD algorithm with linesearch proposed by Malitsky and Pock in 2018 to the setting where the dual space has a Bregman divergence, and the dual function is relatively strongly convex to the Bregman's kernel. This extension yields a new method for OT and WB problems based on smoothing of the objective that also achieves state-of-the-art convergence rates. Finally, we introduce a new Bregman divergence based on a scaled entropy function that makes the algorithm numerically stable and reduces the smoothing, leading to sparse solutions of OT and WB problems. We complement our findings with numerical experiments and comparisons.
\end{abstract}

\maketitle

\begin{small}
  \\[2ex]
\end{small}


\section{Introduction}

We address the problem of computing the Optimal Transport (OT) between two discrete probability distributions and Wasserstein Barycenter (WB) between several (two or more) discrete probability distributions. Optimal Transport is central in machine learning  applications such as classification~\cite{kusner2015word} and unsupervised learning~\cite{arjovsky2017wasserstein,bigot2017geodesic}. On the other hand, the Wasserstein Barycenter problem appears naturally in clustering~\cite{ho2017multilevel} and other problems of imaging~\cite{bonneel2015sliced}.

Discrete Optimal Transport and Wasserstein Barycenter are Linear Programming (LP) problems that can be tackled through interior point methods \cite{lee2015efficient}, or Network Simplex \cite{bonneel2011displacement, flamary2021pot}, which are very efficient (\textit{cf} our experiments in Section~\ref{sec:experiments}). Nevertheless, the increasing size of the data leads to prohibitive large instances for exact solvers. In addition, some common nonlinear extensions might not be tractable by linear optimization techniques. Thus, finding efficient methods to solve these problems by nonlinear optimization is an exciting task in either optimization or computing theory.

Beyond LP, a classical approach to tackle OT is the well-known Sinkhorn algorithm. This method uses an alternating minimization procedure over the dual variables of an entropic regularized version of the OT problem~\cite{cuturi2013sinkhorn,PeyreCuturibook}. Recent papers have analyzed the complexity of Sinkhorn and its variants in terms of the number of arithmetic operations necessary to reach a $\varepsilon-$approximated solution of OT. Most of these rates depend on the error $\varepsilon$, the size of the discrete measures $n$, and the largest element $\|C\|:=\max_{i,j}C_{i,j}$ of the matrix of transportation costs $C\in \mathbb{R}^{n\times n}_+$ (assuming for simplicity, and without loss of generality, $C_{i,j}\ge 0$). In~\cite{altschuler2017near} authors proved the rate $\tilde{O}(\frac{n^2\|C\|^3}{\varepsilon^3})$ for Sinkhorn and its greedy variant named Greenkhorn. Then, in \cite{dvurechensky2018computational} the authors improved the analysis of Sinkhorn and proved the rate $\tilde{O}(\frac{n^2\|C\|^2}{\varepsilon^2})$. The same rate was obtained for Greenkhorn in \cite{lin2019efficient}. Using randomized coordinate descent~\cite{nesterov2012efficiency}, the authors in~\cite{lin2020fixed} proposed the Randkhorn algorithm, a randomized version of Sinkhorn, and prove the rate $\tilde{O}(\frac{n^{7/3}\|C\|^{4/3}}{\varepsilon})$. 

Recently, there has been a growing interest in adapting accelerated convex optimization methods for the entropic regularized OT. Interestingly, these methods achieve better theoretical convergence rates than Sinkhorn-like methods. To the best of our knowledge, the first of these methods was proposed in~\cite{dvurechensky2018computational}. The
authors designed and analyzed the APDAGD method, a linesearch version of the Accelerated Gradient Descent (AGD) scheme~\cite{nesterov2005smooth,tseng2009accelerated}.
Shortly after, Lin et. al.~\cite{lin2019efficient} refined the analysis of APDAGD and prove the rate $\tilde{O}(\frac{n^{5/2}\|C\|}{\varepsilon})$. In the same paper, authors proposed the APDAMD method that slightly generalizes the APDAGD and prove the rate $\tilde{O}(\frac{n^2\sqrt{{\kappa}}\|C\|}{\varepsilon})$, where ${\kappa}=O(n)$ is a factor associated to the strong convexity of the Bregman prox-function with respect to the uniform norm. The same rate was obtained in~\cite{guminov2021combination} using an algorithm based on Accelerated Alternating Minimization (AAM), and also in~\cite{guo2020fast} using Accelerated Primal-Dual
Coordinate Descent.

A different approach was introduced by Jambulapati \textit{et.~al.} in~\cite{jambulapati1906direct}. Here, the authors proposed a new method based on dual-extrapolation~\cite{nesterov2007dual}
and area convexity~\cite{sherman2017area}, and they prove the rate $\tilde{O}(\frac{n^2\|C\|}{\varepsilon})$ (see also \cite{allen2017much}, \cite{cohen2017matrix} and \cite{blanchet2018towards} for similar rates). Despite the theoretical optimality of the rate $O(\frac{n^2\|C\|}{\varepsilon})$ (see \cite{blanchet2018towards}), methods as AGD or AMM outperform in practice the area-convexity approach as is shown by the numerical experiments in~\cite{guminov2021combination} and~\cite{dvinskikh2021improved}. 

The study of the Wasserstein Barycenter problems follows a similar path to OT. For a problem with $m$ marginals, the Iterative Bregman Projection (IBP) method proposed in~\cite{benamou2015iterative} was proven to achieve the rate $\tilde{O}(\frac{mn^2\|C\|^2}{\varepsilon^2})$ in~\cite{kroshnin2019complexity}. Applications of AGD and AAM to WB were also studied in~\cite{kroshnin2019complexity} and~\cite{guminov2021combination}, respectively. Here, authors proved the rate $\tilde{O}(\frac{mn^{5/2}\|C\|}{\varepsilon})$, analogous to the one obtained for OT with the same algorithms. FastIBP algorithm was proposed  in~\cite{lin2020fixed} and authors proved the rate $\tilde{O}(\frac{mn^{7/3}\|C\|^{4/3}}{\varepsilon^{4/3}})$. Recently, in~\cite{dvinskikh2021improved} the idea of area-convexity was extended to WB problems with the near-optimal rate $\tilde{O}(\frac{mn^2\|C\|}{\varepsilon})$. In the same work, authors presented an application of Mirror-Prox~\cite{nemirovski2004prox}) and proved the rate $\tilde{O}(\frac{mn^{5/2}\|C\|}{\varepsilon})$ for the worst-case duality-gap of a saddle-point formulation of WB. Moreover, they provided numerical evidence where Mirror-Prox outperforms the area-convexity method.


Close to Mirror-Prox is the Hybrid Primal-Dual (HPD). This iterative primal-dual method aims to solve a very general class of saddle-points problems and requires half computations per iteration compared to Mirror-Prox. The method was first analyzed in~\cite{chambolle2011first}. Under mild convexity assumptions, the authors proved the rate $O(1/N)$ for the duality-gap after $N$ iterations, and the accelerated rate $O(1/N^2)$ in the presence of strong convexity. Later,
it was revisited in~\cite{chambolle2016ergodic} with an improvement of the analysis and an extension
to general Bregman divergences.

Malitsky and Pock~\cite{malitsky2018first} incorporated a linesearch procedure to the HPD in order to locally estimate the parameter of smoothness at each iteration and allow larger step-sizes. The method with linesearch enjoys the same rates $O(1/N)$ and $O(1/N^2)$ as classic HPD, although the original version is restricted to the Euclidean setting. Recently, in~\cite{jiang2021bregman} it was observed that the analysis of the rate $O(1/N)$ extends to the case in which the dual space is endowed with a Bregman distance while the primal space remains Euclidean. This paper takes one step forward by studying the extension of the rate $O(1/N^2)$. Table~\ref{tab:HPD} brief the results on HPD listed so far.

\begin{table}[ht]
\begin{center}
\begin{tabular}{ccccc}
\hline \hline
Paper & primal \& dual settings & primal \& dual functions & linesearch & rate \\ \hline
\cite{chambolle2011first} & both Euclidean & both convex  & \xmark & $O(1/N)$ \\ 
\cite{chambolle2016ergodic} & both Bregman& both convex & \xmark & $O(1/N)$ \\ 
\cite{chambolle2016ergodic} & both Bregman & convex \& strong convex & \xmark & $O(1/N^2)$ \\ 
\cite{malitsky2018first} & both Euclidean & both convex & \cmark & $O(1/N)$ \\ 
\cite{malitsky2018first} & both Euclidean & convex \& strong convex & \cmark & $O(1/N^2)$ \\ 
\cite{jiang2021bregman} & Euclidean \& Bregman & both convex & \cmark & $O(1/N)$ \\ 
{\bf This paper} & Euclidean \& Bregman & convex \& strong convex & \cmark & $O(1/N^2)$ \\ \hline \hline
\end{tabular}
\caption{Convergence rates for the duality-gap with HPD}
\label{tab:HPD}
\end{center}
\end{table}

\subsection{Contributions}
The contribution of this paper is threefold:
First, we extend the $O(1/N^2)$ result of \cite{malitsky2018first} to the setting where the dual space has a Bregman distance, and the dual function is strongly convex relative to the kernel of the Bregman distance. To do this, we follow closely the proof of Malitsky and Pock by replacing some of the Hilbertian identities with a well-known 3-points inequality that holds in Bregman settings proposed in~\cite{chen1993convergence}. At the end of the proof, we propose a new initialization for the algorithm's parameters that simplifies the analysis and corrects an issue overlooked in the original proof. 

Second, we show how to specialize the HPD into the setting of discrete Optimal Transport and Wasserstein Barycenter problems. To do this, we write these problems in a saddle-point formulation and find tight bounds for the norm of the dual variables.
Interestingly, we find the rates $\tilde{O}(\frac{n^{5/2}\|C\|}{\varepsilon})$ for OT, and $\tilde{O}(\frac{mn^{5/2}\|C\|}{\varepsilon})$ for WB, matching the rates of current state-of-the-art methods in the literature. 

Finally, in Section~\ref{sec:sek} we propose a new scaled entropy kernel that can be accomplished with HPD and other methods in the literature, such as AGD. This new entropy function makes the computation numerically stable, as it avoids computing the logarithm of very small numbers. Moreover, the solutions obtained using the new entropy have smaller support than traditional entropic regularized solutions, which may be of interest for specific applications. We end our paper with a remark (Sec.~\ref{sec:nlp}) on {\it easy} nonlinear extensions of our algorithm and an experimental section comparing our methods with various other techniques from the literature (Sec.~\ref{sec:experiments}).

\section{The Hybrid Primal Dual method} 
We consider the saddle-point problem 
\begin{equation}
\min_{x\in X}\max_{y\in Y} \mathcal{L}(x,y) := g(x)+\la Kx,y\ra-h^*(y),
\label{saddlepoint}
\end{equation}
where $X,Y$ are convex sets of general normed spaces $(\mathcal{X},\|\cdot\|_\mathcal{X})$ (the primal) and $(\mathcal{Y},\|\cdot\|_\mathcal{Y})$ (the dual), respectively. Let $K:\mathcal{X}\mapsto \mathcal{Y}^*$ a linear map such that 
$$L := \sup_{\|x\|_\mathcal{X}\le 1,\|y\|_\mathcal{Y}\le 1} \la Kx,y\ra <\infty,$$ 
and $g,h^*$ proper, lower semicontinuous, convex functions over $X$ and $Y$, respectively. Each space is endowed with a Bregman divergence $D(z,\bar z) = \xi(z)-\xi(\bar z)-\la \nabla \xi(\bar z),z-\bar z\ra,$ 
generated by a kernel or prox-function $\xi$ that is 1-strongly convex with respect to the norm of the space, differentiable in $\textup{int}(\textup{dom}\,\xi)$, and $\|\nabla \xi (x)\|\rightarrow \infty$ when
$x\to\partial(\textup{dom}\,\xi)$
(that is, a strongly convex {\it Legendre function}, as introduced in~\cite{Rockafellar} (Sec.~26), \cite{chen1993convergence}, \cite{bauschkeetal2001}).
We denote by $\xi_\mathcal{X}$ and $\xi_\mathcal{Y}$ the kernel functions on the primal and dual, respectively. Similarly, we denote by $D_\mathcal{X}$ and $D_\mathcal{Y}$ the Bregman divergence, and we assume that
  $X\subset \overline{\textup{dom}}\,\xi_{\mathcal{X}}$ and
  $Y\subset \overline{\textup{dom}}\,\xi_{\mathcal{Y}}$. In addition, we replace without loss of generality $g$ with $g+\delta_X$
    (where $\delta_X$ is the characteristic function of $X$, with value $0$ in $X$
    and $+\infty$ outside)
    and $h^*$ with $h^*+\delta_Y$, and forget in the remaining of the paper about the
    convex constraint sets $X,Y$.

Finally, we define the duality-gap at $(x,y)$ measured at $(\tilde x, \tilde y)$ as $\mathcal{G}_{\tilde{x},\tilde{y}}(x,y):=\mathcal{L}(x,\tilde y)-\mathcal{L}(\tilde x,y)$.

For initial points $x^0,\bar x^0 \in\textup{dom}\,\xi_{\mathcal{X}}, \; y_0\in \textup{dom}\,\xi_{\mathcal{Y}}$, and given non-negative sequences $\{\tau_k\}_k$, $\{\sigma_k\}_k$, $\{\theta_k\}_k$, the main iteration of the HPD method (without linesearch) is given by
\begin{align*}
 y^{k+1} &= \arg\min_{y\in \mathcal{Y}} h^*(y)-\langle K\bar x^k, y\rangle+\frac{1}{\sigma_k}D_{\mathcal{Y}}(y,y^k), \\
 x^{k+1} &= \arg\min_{x\in \mathcal{X}} g(x)+\langle Kx,y^k\rangle+\frac{1}{\tau_k}D_{\mathcal{X}}(x,x^k),\\
 \bar x^{k+1} &= x^{k+1}+\theta_k(x^{k+1}-x^k).
\end{align*}

In the simpler case where $\tau_k\equiv \tau_0$ and $\sigma_k\equiv \sigma_0$ are constant, Theorem~1 in \cite{chambolle2016ergodic} establishes that taking $\theta_k \equiv 1$ and $\tau_0\sigma_0 L^2\le 1$ we have for all $(x,y)\in X\times Y$
\begin{equation}
\mathcal{G}_{x,y}(\hat{x}^N,\hat{y}^N)\le \frac{1}{T_N}\left(\frac{1}{\tau_0}D_\mathcal{X}(x,x^0)+\frac{1}{\sigma_0}D_{\mathcal{Y}}(y,y^0)\right),
\label{simpleHPD}
\end{equation}
where $T_N = N/2$, $\hat{x}^N = \frac{1}{N}\sum_{k=1}^Nx^k$ and $\hat{y}^N = \frac{1}{N}\sum_{k=1}^N y^k$. 

Furthermore, if $g$ is strongly convex relative to the primal prox-function $\xi_\mathcal{X}$ (see Definition 1.2 in \cite{lu2018relatively}), then Theorem~4 (together with footnote~1, p. 268)
in \cite{chambolle2016ergodic} establishes that taking $\theta_{k+1} = \frac{1}{\sqrt{1+\gamma\tau_k}}$, $\tau_{k+1} = \theta_{k+1}\tau_k,$ and  $\sigma_{k+1} = \sigma_k/\theta_{k+1}$, equation \eqref{simpleHPD} holds with $T_N \gtrapprox \gamma N^2/(L^2\sigma_0)$, where the notation $\gtrapprox$ means that the inequality holds up to constant factors and for a sufficiently large $N$. We remark that a symmetric algorithm and result also holds in the case that $h^*$ is strongly convex relative to $\xi_\mathcal{Y}$, which is the case of interest for the next section.

\begin{remark}[Erratum to~\cite{chambolle2016ergodic}]
  In~\cite{chambolle2016ergodic} it is claimed the convergence of the
  iterates holds in general, yet the proof works only whenever the
  prox-function is differentiable on the whole domain (as is the case
  for instance for the scaled entropy of Sec.~\ref{sec:sek}) or when
  the limit points lie in the interior (as in the entropy-regularized
  problem of Sec.~\ref{sec:regOT}). This is observed for instance in
  Theorem~4 in~\cite{lee2021lastiterate}, or~\cite{LeeLastIterate}.
  The algorithm also extends with easy adaption to strongly convex
  prox-functions, possibly non-smooth, but with locally bounded
  subgradients, as pointed out in~\cite{LanZhou18} where a variant of
  the above Bregman primal-dual descent is introduced. This allows to
  consider prox-functions with finite gradients on the boundary of
  their domain.  In these cases the convergence proof
  of~\cite{chambolle2016ergodic}, which uses the fact that limits of
  iterates are also fixed points of the algorithm (which is not
  correct if the prox-function has infinite gradients at such points
  as the algorithm is then undefined), may be adapted.  The
  authors of~\cite{chambolle2016ergodic} apologize for letting this pass through the publication
  process and warmly thank Chung-Wei Lee from U.~South~Carolina for
  pointing out this issue.
\end{remark}

\subsection{Linesearch}
\label{sec:linesearch}

Algorithm~\ref{Alg:PDls} presents a primal-dual method with linesearch slightly different from the original proposed by \cite{malitsky2018first}. The results in~\cite{malitsky2018first} guarantee the rates $O(1/N)$ and $O(1/N^2)$ for the particular  case $D_{\mathcal{X}}(x,\bar x) = \frac{1}{2}\|x-\bar x\|^2_2$ and $D_{\mathcal{Y}}(y,\bar y) = \frac{1}{2}\|y-\bar y\|^2_2$, which is referred to as the Euclidean setting. Recently, \cite{jiang2021bregman} proved that the
$O(1/N)$ algorithm can be extended to the setting where the Bregman distance $D_\mathcal{Y}$ has a non-Euclidean kernel.

\begin{small}
\begin{algorithm}[ht]
\begin{algorithmic}[1]
 \caption{HPD method with linesearch}\label{Alg:PDls}
 \REQUIRE $x^1 =x^0\in \textup{dom}\,\xi_{\mathcal{X}}, \; y^1\in \textup{dom}\,\xi_{\mathcal{Y}}$, $\beta_0>0$, 
 $\tau_0=\frac{1}{\sqrt{\beta_0}L}$, $\gamma\ge 0$, 
$\theta_0>0$,
$\rho\in (0,1)$, set $k=1$.  
\REPEAT[mean-loop]
\STATE 
 $\tau_k \leftarrow \tau_{k-1}\sqrt{1+\theta_{k-1}}/\rho, \quad \beta_{k} \leftarrow \frac{\beta_{k-1}}{1+\gamma\beta_{k-1}\tau_{k-1}}$ \alglabel{step:taubeta}
\REPEAT[inner-loop: linesearch]
\STATE \alglabel{linesearch} $\tau_k \leftarrow \rho\tau_k, \quad \theta_k \leftarrow \tau_k/ \tau_{k-1}, \quad \sigma_k \leftarrow \beta_k\tau_k$
\STATE $\bar x^k \leftarrow x^k+\theta_k(x^k-x^{k-1})$ 
\STATE $y^{k+1} \leftarrow  \displaystyle{\arg\min_{y\in \mathcal{Y}}}\; h^*(y)-\langle K\bar x^k,y\rangle +\frac{D_{\mathcal{Y}}(y,y^{k})}{\sigma_k}$ \alglabel{dualprox}
\STATE $x^{k+1} \leftarrow \displaystyle{\arg\min_{x\in \mathcal{X}}}\; g(x)+\langle Kx,y^{k+1}\rangle +\frac{\|x-x^{k}\|^2_2}{2\tau_{k}}$ \alglabel{primalprox}
\UNTIL{\alglabel{MPrule} $\frac{1}{2}\|x^{k+1}-\bar x^k\|^2_2+\frac{1}{\beta_k}D_{\mathcal{Y}}(y^{k+1},y^k)+\tau_k\langle  K(x^{k+1}-\bar x^k),y^{k+1}-y^k\rangle\ge 0$ }
\STATE $k \leftarrow k+1$
\UNTIL{stopping criteria}
\end{algorithmic}
\end{algorithm}
\end{small}

The following Theorem extends the results of \cite{malitsky2018first} and \cite{jiang2021bregman} by observing that if $h^*$ is strongly convex relative to $\xi_\mathcal{Y}$, then the accelerated rate $O(1/N^2)$ can be achieved.
Its proof is given in~Appendix~\ref{sec:proofMP}.
\begin{theorem}\label{thm:MP}
Let $(x^k,\bar x^k, y^k,\tau_k,\sigma_k,\theta_k)_k$ the sequences generated by Algorithm~\ref{Alg:PDls}, and $T_N = \sum_{k=1}^N\tau_k$. For any $(x,y)\in \mathcal{X}\times \textup{dom}\, \xi_{\mathcal{Y}}$ it holds
\begin{equation}
\mathcal{G}_{x,y}(\hat x^N,\hat y^N) \le \frac{1}{T_N}\bigg(\frac{1}{2}\|x-x^1\|^2_2+  \frac{1}{\beta_1}D_{\mathcal{Y}}(y,y^1)-\tau_1\theta_1\big(\mathcal{L}(\hat x^N,y)-\mathcal{L}(x^0,y)\big)\bigg),
\label{AMP}
\end{equation}
where the ergodic outputs are now $\hat x^N = \frac{\tau_1\theta_1x^0+\sum_{k=1}^N\tau_k\bar x^k}{\tau_1\theta_1+T_N}$ and $\hat y^N = \sum_{k=1}^N \frac{\tau_k}{T_N}y^{k+1}$.

The sequence $\beta_k\equiv \beta$ with $\tau_0\sqrt{\beta}L=1$ yields $T_N \gtrapprox N$. 
 If $h^*$ is $\gamma$-strongly convex relative to $\xi_{\mathcal{Y}}$
  and $0< \gamma \le L/\rho$,
  then the choice $\beta_0>0$, $\tau_0=1/(\sqrt{\beta_0}L)$ and $\theta_0=\gamma\sqrt{\beta_0}/L$
  yields $T_N \gtrapprox \gamma N^2 / L^2$.

If in addition the Lagrangian $\mathcal{L}(x,y)$ is linear with respect
to feasible points $x$, then one can define
$\hat x^N = \sum_{k=1}^N\frac{\tau_k}{T_N}\bar x^k$ and~\eqref{AMP} simplifies to
\begin{equation}
\mathcal{G}_{x,y}(\hat x^N,\hat y^N) \le 
\frac{1}{T_N}\bigg(\frac{1}{2}\|x-x^1\|^2_2+  \frac{1}{\beta_1}D_{\mathcal{Y}}(y,y^1)\bigg).
\label{AMPlin}
\end{equation}
\end{theorem}

\begin{remark}[Difference with respect to~\cite{malitsky2018first}] The main differences with respect to the original algorithm are the setup for the strongly convex setting and the linesearch procedure (steps~\ref{linesearch} to~\ref{MPrule} in Algorithm~\ref{Alg:PDls}). We choose $\theta_0 = \gamma \sqrt{\beta_0}/L$ instead of $\theta_0 = 1$,
    and more importantly, we fix the initial value of $\tau_k$ before the linesearch as $\tau_k = \tau_{k-1}\sqrt{1+\theta_{k-1}}$, while the original algorithm proposes to choose $\tau_k \le \tau_{k-1}\sqrt{1+\theta_{k-1}}$.
    These modifications seem necessary in our analysis (see Appendix~\ref{sec:proofMP})
    to prove the relation $\tau_k \ge \frac{1}{L\sqrt{\beta_k}}$ for all $k\ge 1$.
    This property was used without a proof in~\cite{malitsky2018first}, yet the authors
    seem to have overlooked that it is unclear in the strongly convex case.

In the linesearch, the stopping criteria for the inner-loop (step~\ref{MPrule} in Algorithm~\ref{Alg:PDls}) was first proposed in \cite{jiang2021bregman}. This new condition yields larger steps than the original one, yet it requires to compute $x^{k+1}$ at each iteration of the inner-loop. Hence, it is applicable only when the computation of $x^{k+1}$ is relatively cheap, which will be the case in  our implementations. In other cases one can revert to the termination criterion of~\cite{malitsky2018first}. Finally, we prove in the Appendix~\ref{termination} that the linesearch always terminates, and the overall number of inner-loop iterations of after $N$ main-loop iterations is $O(N/|\ln \rho|)$, therefore does not affect the rate of convergence. 
\end{remark}

\section{Application to Optimal Transport and Wasserstein Barycenter problems}
\label{Sec:Specialization}

We start this section by formally defining the discrete Optimal Transport problem.

Let $C\in \mathbb{R}^{n\times n}_+$ and $\mu,\nu\in \mathbb{R}^n_+$ two discrete probability measures, i.e. $\onen^\top \mu = \onen^\top \nu = 1$ where $\onen = (1,1,...\,,1)^\top\in\mathbb{R}^n$. We consider the discrete Optimal Transport problem 
\begin{equation}
\min_{X\in \mathbb{R}^{n\times n}_+}\left\{\langle C,X\rangle \; :\;X\onen=\mu, X^\top \onen = \nu\right\}.
\label{OT}
\end{equation}
The saddle-point formulation of \eqref{OT} is
\begin{equation}
\min_{X\in {\Delta}}\max_{u,v}\la C,X \ra + \la u ,\mu-X\onen \ra+\langle v,\nu-X^\top\onen\rangle,
\label{OTsaddle}
\end{equation}
where ${\Delta} = \{X\in \mathbb{R}^{n\times n}_+: \onen^\top X \onen = 1\}$ denotes the Simplex of square matrices of size $n$.

Next, we would like to bound the dual variables as $\|(u,v)\|\le \lambda$, where $\|\cdot\|$ stands for the norm of the maximum, and $\lambda$ is a positive constant. Lemma~2 in \cite{jambulapati1906direct} guarantees that $\lambda = 2\|C\|$ is large enough to bound the dual variables. Here, we improve this estimate by determining the tightness value of the constant $\lambda$. See Appendix~\ref{sec:halfC} for the proof.
\begin{lemma}\label{thm:halfC}
Assume $\min_{i,j} C_{i,j}=0$. Then problem \eqref{OTsaddle} admits a solution $(u^*,v^*)$ with $\|(u^*,v^*)\|\le \frac{\|C\|}{2}$. Moreover, there exist a matrix $C$ and marginals $\mu,\nu$ such that for any dual solution $(u^*,v^*)$ of problem \eqref{OTsaddle} it holds $\|(u^*,v^*)\|\ge \frac{\|C\|}{2}$. 
\end{lemma} 

In equations \eqref{saddlepoint} and \eqref{OTsaddle} the role of primal and dual is interchangeable. As the dual space $\mathcal{Y}$, we consider the space of square matrices of size $n$ with the norm $\|X\|_1 = \sum_{i,j} |X_{ij}|$. In this setting it is classical to use the entropy kernel given by $\xi_{\mathcal{Y}}(X) = \sum_{ij}X_{ij}\ln X_{ij}$, that defines the Bregman distance 
$$D_{\mathcal{Y}}(X,\bar X) = \langle X,\ln X-\ln \bar X\rangle.$$
We recall that Pinsker's inequality (see for instance~\cite{Tsybakov}) expresses that the entropy  $\xi_\mathcal{Y}(X)$ is 1-strongly convex with respect to the norm $\|\cdot\|_1$ as long as $X$ remains in the simplex, and the associated Bregman distance it is nothing but the Kullback-Leibler divergence. For the primal space $\mathcal{X}$, we consider $\mathbb{R}^n\times \mathbb{R}^n$ with the usual Euclidean norm denoted by $\|\cdot\|_2$ and the kernel $\xi_\mathcal{X}(u,v) = \frac{1}{2}\|(u,v)\|^2_2$. 

Comparing the structure of \eqref{saddlepoint} and \eqref{OTsaddle}, we observe that
$g(u,v) =  -\langle (u,v),(\mu,\nu)\rangle$, which is linear;  $h^*(X) =\langle C,X\rangle + {\delta}_\Delta(X)$, which is the sum of a linear and a characteristic function; and $K(u,v) = u\otimes \onen+\onen\otimes v$, where $\otimes$ denotes the tensor product. Moreover, we can compute explicitly
$$L= \sup_{\|X\|_1 \le 1, \;\; \|(u,v)\|_2 \le 1} \langle X, u\otimes \onen+\onen\otimes v \rangle  = \sqrt{2}.$$

\subsection{Non-regularized optimal transport} \label{sec:non-regOT}

We first analyze the non-regularized version of the Optimal Transport problem. Consider problems in steps \ref{dualprox} and \ref{primalprox} of Algorithm~\ref{Alg:PDls} which  are called the dual and primal proximal problems, respectively. For the saddle-point formulation of OT \eqref{OTsaddle}, these problems admits explicit solutions, namely 
\begin{equation}
\begin{aligned}
X^{k+1} &= \mathcal{N}\left(X^k\exp \left(-\sigma_k(C - \bar u^k \otimes\onen-\onen\otimes \bar v^k)\right)\right), \\
u^{k+1} &= \Pi_{[-\lambda,\lambda]}\left(u^k+\tau_k(\mu - X^{k+1}\onen)\right), \\
v^{k+1} &= \Pi_{[-\lambda,\lambda]}\left(v^k+\tau_k(\nu - (X^{k+1})^\top\onen)\right), \end{aligned}
\label{solutionsOT}
\end{equation}
where $\mathcal{N}(X) = X/\|X\|_1$ denotes the renormalization of matrix $X$, and $\Pi_{[-\lambda,\lambda]}$ denotes the component-wise projection on the interval $[-\lambda,\lambda]$. We remark that the computation of $(u^{k+1},v^{k+1})$ is practically inexpensive.

As the starting points of Algorithm~\ref{Alg:PDls} we take $X^1 = \frac{1}{n^2}(\onen\otimes \onen)$, and $u^0=v^0=u^1=v^1=0_{\mathbb{R}^n}$. With this, it is easy to check that $D_{\mathcal{Y}}(X,X^1)  \le 2\ln n$ for all $X\in {\Delta}$, and $\frac{1}{2}\|(u,v)-(u^0,v^0)\|_2^2\le n\lambda^2$ for all pairs $(u,v)$ such that $\|(u,v)\|\le \lambda$. Note that alternatively we may choose $X^1 = \mu\otimes \nu$ which yields similar bounds. 

  For slightly improved convergence rates, we rely on
  the following Lemma which improves the estimate of Lemma 7 in \cite{altschuler2017near} by a factor 2. See Appendix~\ref{sec:approx2}
  for its proof.
\begin{lemma}\label{thm:approx2}
  Algorithm~2 in~\cite{altschuler2017near} (``\textsc{Round}''),
  denoted $\mathcal{A}$, when receiving a matrix $X\in {\Delta}$
  and marginals $\mu,\nu$, returns in $O(n^2)$ arithmetic operations a matrix $\mathcal{A}(X)\in{\Delta}$ such that $\mathcal{A}(X)\onen = \mu$, $\mathcal{A}(X)^\top\onen = \nu$, and
$$\|X-\mathcal{A}(X)\|_1\le \|\mu-X\onen\|_1+\|\nu - X^\top\onen\|_1.$$
\label{alts}
\end{lemma}
\begin{proposition}\label{thm:boundMP}
  The HPD algorithm with linesearch finds an $\varepsilon-$optimal solution to \eqref{OT} in $\tilde{O}(\frac{n^{5/2}\|C\|}{\varepsilon})$ arithmetic operations.
\end{proposition}
\begin{proof} Let $\hat X^N$ and $(\hat u^N,\hat v^N)$ the solutions reported by the algorithm after $N$ iterations. As the functions involved in OT are linear we have $\gamma = 0$. 

Denote $X^*$ as the true optimum of~\eqref{OT}. Using Cauchy-Schwarz's inequality and Lemma~\ref{alts}, we have
\begin{align*}
\langle C,\mathcal{A}(\hat X^N)-X^*\rangle & = \langle C,\mathcal{A}(\hat X^N)-\hat X^N+\hat X^N-X^*\rangle \\ 
& \le \langle C,\hat X^N-X^*\rangle +\lambda(\|\mu-\hat X^N\onen\|_1+\|\nu - (\hat X^N)^\top\onen\|_1).
\end{align*}
The last right hand side is exactly the duality-gap at $X = X^*$, $u = \lambda\mbox{sign}(\mu-\hat X^N\onen)$, and $v = \lambda\mbox{sign}(\nu-(\hat X^N)^\top\onen))$. Then, using the bound of Theorem~\ref{thm:MP}, and optimizing on $\beta,\tau_0$ with the constraint $\tau_0\sqrt{\beta} L = 1$ leads to $\tau_0 = \frac{\lambda\sqrt{n}}{L\sqrt{2\ln n}}$ and $\beta = \frac{2\ln n}{\lambda^2n}$. Replacing and recalling that $T_N \gtrapprox N$, we obtain the rate $\tilde{O}(\frac{\sqrt{n}\|C\|}{\varepsilon})$ for the number of iterations. The result follows by noticing that each iteration requires component-wise operations over a square matrix of size $n$, then $O(n^2)$ arithmetic operations per iteration.
\end{proof}

\subsection{Entropy regularized Optimal Transport}
\label{sec:regOT}

In this section we consider the entropic regularized version of the saddle-point formulation of OT. To do this, we add an entropy regularization term to~\eqref{OTsaddle}, so now $h^*(X) = \la C,X\ra+\gamma\la X,\ln X\ra + \delta_\Delta(X)$, with $\gamma>0$. Clearly the new function $h^*$ is $\gamma-$strongly convex on ${\Delta}$ relative to the entropy kernel and we obtain the entropic regularized saddle-point formulation for OT
\begin{equation}
\min_{X\in {\Delta}}\max_{u,v}\la C+\gamma\ln X,X \ra + \la u ,\mu-X\onen \ra+\langle v,\nu-X^\top\onen\rangle.
\label{OTsaddle-reg}
\end{equation}
 Note that the same formulation was recently used in \cite[Sec.~5.2]{silveti2021stochastic} to tackle inverse problems on discrete measures.

We emphasize that the new entropy term in equation \eqref{OTsaddle-reg} only affects the solution of the proximal problem associated to the matrix $X$, which is now given by 
$$X^{k+1} = \mathcal{N}\left(\exp\left(\frac{1}{1+\sigma_k\gamma}\left(\ln X^k - \sigma_k(C-\bar u^k\otimes \onen-\onen\otimes \bar v^k)\right)\right)\right).$$ 

\begin{proposition}\label{thm:boundaccel}
Let $\varepsilon>0$ and $\gamma = \frac{\varepsilon}{4\ln n}$. The HPD algorithm for $\gamma$ entropy regularized OT requires $O(\frac{n^{5/2}\|C\|}{\varepsilon})$ arithmetic operations to achieve a $\varepsilon-$optimal solution of \eqref{OT}. 
\end{proposition}
\begin{proof}
  In order to bound the number of iterations we proceed similarly to
  the proof of~\cite[Thm~4]{dvurechensky2018computational}. Let $X^*_\gamma$ the optimal solution of the entropy regularized problem, then
\begin{align*}
\langle C ,\mathcal{A}(\hat X^N)&-X^*\rangle  = \langle C,\mathcal{A}(\hat X^N)-\hat X^N\rangle+\langle C,\hat X^N-X^*_\gamma\rangle+\langle C,X^*_\gamma-X^*\rangle \\
 \le\ & \langle C,\hat X^N-X^*_\gamma\rangle+\gamma(\xi_{\mathcal{Y}}(\hat X^N)-\xi_{\mathcal{Y}}(X^*_\gamma)) 
+\lambda(\|\mu-\hat X^N\onen\|_1+\|\nu - (\hat X^N)^\top\onen\|_1)
  \\  & 
 +\langle C,X^*_\gamma-X^*\rangle +\gamma(\xi_{\mathcal{Y}}(X^*_\gamma) 
-\xi_{\mathcal{Y}}(X^*))+\gamma(\xi_{\mathcal{Y}}(X^*)-\xi_{\mathcal{Y}}(\hat X^N)).
\end{align*}
Here, the inequality is due to Lemma \ref{alts}. In the last expression, the first line correspond to the duality-gap at $X = X^*_\gamma$, and $(u,v)=(\lambda\mbox{sign}(\mu-\hat X^N\onen),\lambda\mbox{sign}(\nu-(\hat X^N)^\top\onen))$. In the second line, the  sum of the first two terms is $\le 0$ because of the optimallity of $X^*_\gamma$, while the third term is $\le 2\ln n$ cause $\xi_\mathcal{Y}(X) \in [-2\ln n, 0]$ for all $X\in{\Delta}$. We recall that for OT the function
$\mathcal{L}(\cdot,y)$ is linear so the bound in equation \eqref{AMPlin} implies
\begin{equation}
\langle C, \mathcal{A}(\hat{X}^N)-X^*\rangle\le \frac{1}{T_N}\left(n\lambda^2 + \frac{1}{\beta_1}2\ln n\right)+2\gamma \ln n.
\label{AMP2-OT}
\end{equation}
Given a tolerance $\varepsilon>0$, we take $\gamma = \frac{\varepsilon}{4\ln n}$. On the other hand, we can choose a large $\beta_1$ in order to minimize the factor in parenthesis, or we may choose $\beta_1 = O(\frac{2\ln n}{n\lambda^2})$ in order to equate the two terms in the sum. Note that smaller values of $\beta_1$ allow larger initial stepsizes $\tau_1$ as the initialization implies $\tau_1 = 1/\sqrt{2\beta_1}$ (recalling $L=\sqrt{2}$).
We discuss the impact of different choices for $\beta_1$ in the experimental Section~\ref{sec:experiments}. Finally, using Theorem~\ref{thm:MP} we get $T_N \gtrapprox \gamma \rho^2/32 N^2$, so at least $\tilde{O}(\frac{\sqrt{n}\|C\|}{\rho\varepsilon})$ iterations are required to reach the precision $\varepsilon$.
\end{proof}

\begin{remark}
Interestingly, the entropy regularized algorithm's convergence rate is the same as for the method without regularization. In the experimental section, we observe that the regularized version performs better. On the other hand, the non-regularized algorithm does not depend on the prefixed tolerance $\varepsilon$, so it might be helpful when the tolerance is not defined beforehand or whenever an early stopping criterion is required.
\end{remark}

\subsection{Application to Wasserstein Barycenter}
\label{Sec:SpecializationWB}

Given discrete distributions $\mu_1,...\,,\mu_m$,  cost matrices $C_1,...\,,C_l$, and weights $w_1,...\,,w_m\ge 0$ with $\sum_i w_i=1$, the Wasserstein Barycenter problem is to determine a central distribution $\nu$ that minimize the weighted OT distance to the given distributions. Similar to OT, we consider the the saddle-point formulation of this problem
\begin{equation}
\min_{\nu\in {\Delta},\;  X^l\in {\Delta}} \max_{u,v} \sum_{l=1}^m w_l\left(\langle C_l,X_l\rangle + \langle u_l,\mu_l-X_l\onen\rangle+\langle v_l,\nu-X_l^\top\onen\rangle\right).
\end{equation}
Here, it may be convenient to eliminate the barycenter (variable $\nu$) by replacing it with one of the constraints, for instance  $\nu = X_m^\top \onen$. By doing this, and defining conveniently $v^m$  such that $\sum_{l=1}^{m} w_lv_l = 0$, we have
$$\sum_{l=1}^m w_l\la v_l, \nu-X_l^\top\onen\ra = -\sum_{l=1}^{m} w_l\la v_l,X_l^\top\onen \ra.$$

Analogously to OT, we choose the dual space to be the space of tensors $X = (X_1,...\,,X_m) \in \Delta^m$ with the weighted Bregman distance given by $D_\mathcal{Y}(X,\bar X) = \sum_l w_l \langle X_l,\ln X_l-\ln \bar X_l\rangle$ (which correspond to the weighted entropy kernel $\xi_{\mathcal{Y}}(X) = \sum_l w_l\langle X_l,\ln X_l\rangle$). Similarly, we consider the weighted Euclidean prox $\xi_{\mathcal{X}}(u,v) = \frac{1}{2}\sum_l w_l\|(u_l,v_l)\|_2^2$. 

Under the considerations above, we recover explicit solutions for the proximal dual and primal problems similar to the ones of equation~\eqref{solutionsOT}, namely, for the dual 
\begin{equation*}
X^{k+1}_l = \mathcal{N}\left(X^k_l\exp \left(-\sigma_k(C_l - \bar u^k_l \otimes\onen-\onen\otimes \bar v^k_l\right)\right), \qquad \mbox{for all }l=1,...,m,
\end{equation*}
and then for the primal,
\begin{equation*}
u^{k+1}_l = \Pi_{[-\lambda,\lambda]}\left(u^k_l+\tau_k(\mu^l - X^{k+1}_l\onen)\right),  \qquad  \mbox{for all } l = 1,...,m,
\end{equation*}
\begin{equation*}
v^{k+1}_l = \Pi_{[-\lambda,\lambda]}\left(v^k_l+\tau_k(X^{k+1}_m - X^{k+1}_l)^\top\onen\right), \qquad \mbox{for all } l=1,...,m-1.
\end{equation*}
Finally, the last variable $v_m$ updates as $v^{k+1}_m = -\frac{1}{w_m}\sum_{l=1}^{m-1}w_lv^{k+1}_l$.

The number of iterations necessary to reach the precision $\varepsilon$ can be bounded in the same way as for OT leading to $\tilde{O}(\frac{\sqrt{n}\|C\|}{\varepsilon})$ iterations. We remark that each iteration requires $m$ independent component-wise operations over a square matrix of size $n$, hence a total of $\tilde{O}(\frac{mn^{5/2}\|C\|}{\varepsilon})$ arithmetic operations, matching the rates of the recent methods in the literature.

\begin{remark}[Fixed marginal] \label{fixmarginal} In the settings described so far we introduce variables $(u,v)$ that correspond to constraints $X\onen = \mu$ and $X^\top \onen = \nu$, respectively. We remark that it is possible to eliminate one of these variables by keeping the respective constraints in the optimization problem. For instance, we eliminate the variable $u$ by changing the constraint $X\in {\Delta}$ for $X\in {\Delta}_{\mu} := \{X\in \mathbb{R}^{n\times n}_+: X\onen = \mu\}$. We refer to this as the ``fixed marginal approach''. In this case, the renormalization for the variable $X$ is now given by $[\mathcal{N}_{\mu}(X)]_{ij} = \mu_i\frac{X_{ij}}{\sum_{j'}X_{i,j'}}$. Note that the reduction on the dimension also has consequences for the computation of the parameter $L$, which goes down from $\sqrt{2}$ to $1$ in OT problems. Although these modifications only improve on the constants, we will see in the experimental section that fixing a marginal leads to smaller values of the duality-gap and improves the overall performance at no additional cost. 
\end{remark}

\section{Extensions and variants}

\subsection{Scaled entropy kernel}\label{sec:sek}

In this section we introduce a new scaled entropy kernel. This new entropy function allows to solve the Optimal Transport problem in a computational stable way by bounding the elements of the matrix $X$ away from 0. Hence, we avoid computing the logarithm of very small numbers, which is an usual drawback of methods based on entropic regularization.
Given ${\delta}\in (0,1)$,  we define the function:
\begin{equation}
\xi_{\delta}(X) = \frac{1}{(1-\delta)^2}\xi\left((1-{\delta})X+\frac{{\delta}}{n^2}\one_{n\times n}\right),
\label{scaledentropy}
\end{equation}
where $\xi(X)=\sum_{ij}X_{ij}\ln X_{ij}$ is the usual entropy kernel. The natural domain for $\xi_{\delta}$ is $\{X\in \mathbb{R}^{n\times n}: X_{ij}\ge -\frac{\delta}{n^2}, \; \sum_{ij} X_{ij} = 1\}$, however, we restrict to matrices with non-negative entries by adding to the problem the constraint
$X\ge 0$. 

The modified entropy inherits the properties of a kernel function. In particular, $\xi_\delta$ is 1-strongly convex w.r.t. to the norm 1. Moreover, denoting $X^{\delta} = (1-\delta)X+\frac{\delta}{n^2}\one_{n\times n}$, we have that the Bregman distance associated to $\xi_\delta$ is $D_\delta(X,\bar X) = \frac{1}{(1-\delta)^2}D(X^\delta,\bar X_\delta)$, where $D$ is the usual Kullback–Leibler divergence generated by $\xi$. Finally, we remark that if $X\ge 0$ is feasible, i.e. $X\onen = \mu$ and $X^\top\onen = \nu$, then $X^{\delta}\ge \frac{\delta}{n^2}$ and satisfies the marginal constraints $X^\delta\onen = \mu^{\delta}$ and $(X^\delta)^\top = \nu^\delta$ for the $\delta-$regularized marginals $\mu^{\delta} = (1-\delta)\mu+\frac{\delta}{n}\onen$ and $\nu^\delta = (1-\delta)\nu+\frac{\delta}{n}\onen$.   

At first sight, the idea of the scaled entropy seems nonsense as the usual entropy is precisely introduced to ensure the positivity constraint, while the new scaled entropy makes the computation of the prox difficult again. However, as we show in the next section, for OT and WB, these difficulties are overcome by an iterative subroutine that computes an approximated solution to the prox in a few iterations. Moreover, we show in the numerical experiments that this approach leads to sparse solutions of OT and WB that may be of interest in applications (see for instance~\cite{rabin2014adaptive,essid2018quadratically,blondel2018smooth}).

\subsubsection{Optimal Transport with scaled entropy}\label{sec:fixed}

We start by recalling the dual proximal problem (step \ref{dualprox} in Algorithm \ref{Alg:PDls}). By incorporating the scaled entropy kernel, and getting rid of the index $k$, solving this prox problem is equivalent to solve
\begin{equation}
\min_{X^{\delta} \in {\Delta}{\delta}} \la C-\bar{u}\otimes\onen-\onen\otimes \bar v, X^{\delta}\ra+\frac{1}{\sigma}D(X^{\delta},\bar X^{\delta}),
\label{scaledprox}
\end{equation}
where $\bar X^{\delta}$ is the scaled transformation of the last iterate, and ${\Delta}{\delta} = \{X\in {\Delta}: X_{ij}\ge \frac{{\delta}}{n^2}\}$. 

These new constraints guarantee that the components of $X^{\delta}$ remain bounded away from 0, nevertheless, there is no longer an explicit solution for problem~\eqref{scaledprox}. In the following, we
detail a simple and efficient method to compute an accurate solution of~\eqref{scaledprox}.

Consider the optimality conditions of \eqref{scaledprox}: there exist multipliers $\varrho\in \mathbb{R}$, $\varrho_{ij}\ge 0$ such that:
\begin{equation}
X^{\delta}_{ij} = \bar X_{ij}^{\delta} \exp\left( -\sigma(C_{i,j}-\bar{u}_i-\bar{v}_j)-\varrho+\varrho_{ij}\right).
\label{scaledOTcond}
\end{equation}
Furthermore, the complementary conditions implies that if $\varrho_{ij}>0$ then $X^{\delta}_{ij} = \frac{{\delta}}{n^2}$. Denote $s = \exp(\varrho)$, one has therefore
$$X^{\delta}_{ij} = \max\left\{\frac{1}{s}\bar X_{ij}^{\delta} \exp\left( - \sigma(C_{i,j}-\bar{u}_i-\bar{v}_j)\right), \frac{{\delta}}{n^2}\right\}.$$
Summing over $i$ and $j$, and recalling the constraint $\sum_{ij} X^{\delta}_{ij} = 1$ yields
$$s = \sum_{i,j}\max\left\{\bar X_{ij}^{\delta} \exp\left( - \sigma(C_{i,j}-\bar{u}_i-\bar{v}_j)\right),s\frac{{\delta}}{n^2}\right\}.$$
Let us define $Z_{ij} := \bar X_{ij}^{\delta} \exp\left( - \sigma(C_{i,j}-\bar{u}_i-\bar{v}_j)\right)$. Then the problem reduces to computing the value $s>0$, which is the unique root of the piecewise affine and concave real function 
\begin{equation}
F(s):= s-\sum_{ij} \max\left\{Z_{ij}, s\frac{{\delta}}{n^2}\right\}.
\label{Newtonfunction}
\end{equation}
Moreover, one has $F(0) = -\sum_{ij}Z_{ij}<0$ and $F(n^2\max_{ij}Z_{ij}/{\delta}) = (1-{\delta})s\ge 0$, hence $s$ can be approximated using Newton's method starting from $s^0=0$. We remark that $F'(s)$ has the simple expression
$$F'(s) = 1-\frac{{\delta}}{n^2}\left\vert \left\{ (i,j): Z_{ij}\le s^k\frac{{\delta}}{n^2}\right\}\right\vert,$$
where $|\cdot|$ denotes the cardinality of a set. Then, the iteration of Newton's method reads:
$$s^{k+1} = \frac{\sum_{(i,j): Z_{ij}>s^k{\delta}/n^2}Z_{ij}}{1-{\delta}+\frac{{\delta}}{n^2}\vert\{(i,j): Z_{ij}>s^k\frac{{\delta}}{n^2}\}\vert}.$$

Alternatively, the value of $s$ can be approximated using the classical Bamach-Picard iteration to solve the fixed point problem $s = Ts$ for the operator $Ts = \sum_{ij} \max\left\{Z_{ij}, s\frac{{\delta}}{n^2}\right\}$.  Note that for all $s,s'\ge 0$ we have that $|Ts'-Ts| \le {\delta}|s-s'|$, so $T$ is a ${\delta}$-contraction and we have the estimate $|s^k-s^*|\le {\delta}^k|s^*-s^0|$, where $s^*$ is the unique fixed point of $T$. For instance, with ${\delta} = 0.1$, we reach the tolerance $\epsilon$ after $\left\lceil \log_{10} \frac{|s^*|}{\epsilon} \right\rceil$ iterations, so that in this case the number of significant digits equals the number of iterations. In
practice, the Newton method achieves a good precision in even fewer iterations.

\begin{remark}
  One can also be tempted to solve the OT problem by alternating maximization on
  the dual of the scaled entropy-regularized assignment problem, in a sort
  of ``scaled Sinkhorn algorithm''. This works well in
  practice,
  but a (good) theoretical bound on the number of iterations needed for solving
  each maximization is missing.
\end{remark}

\subsubsection{Accelerated Gradient Descent with scaled entropy}\label{sec:appscaled}

The scaled entropy can also be incorporated into other methods in the literature which rely on entropy smoothing,
as is the case for the Accelerated Gradient Descent.

Let us consider the Optimal Transport problem:
$$\min_{X\ge 0}\{\langle C,X\rangle: X\onen = \mu, X^\top\onen=\nu\}.$$ 
Now, we change the variable $X$ for $X^{\delta} = (1-{\delta})X+\frac{{\delta}}{n^2}$. Straightforward computations show that the new problem is equivalent (in the sense that they have the same optimum) to
$$\min\left\{\langle C,X^{\delta}\rangle: X^{\delta}\ge \frac{{\delta}}{n^2}, X^{\delta}\onen = \mu^{\delta}, (X^{\delta})^\top\onen=\nu^{\delta}\right\}.$$ 

Now we add the entropic regularization term $\gamma\la X^{\delta},\ln X^{\delta}\ra$, and, following \cite{dvurechensky2018computational} find the dual:
\[\max_{u,v}\;  \varphi(u,v):=\langle u,\mu^{\delta}\rangle+\langle v,\nu^{\delta}\rangle + V(u,v),
\makebox[0cm][l]{\quad with:}\]
\begin{equation}
V(u,v):= \min_{X^{\delta}\in {\Delta}{\delta}}\langle C-u\otimes \onen-\onen\otimes v,X^{\delta}\rangle + \gamma\langle X^{{\delta}},\ln X^{\delta}\rangle.
\label{AGDdeltaprox}
\end{equation}

If we denote by $X^{\delta}(u,v)$ the solution of problem \eqref{AGDdeltaprox}, then we have $\nabla_u \varphi(u,v) = \mu^{\delta}-X^{\delta}(u,v)\onen$, and $\nabla_v \varphi(u,v) = \nu^{\delta}-X^{\delta}(u,v)^\top\onen$.

The main difficult is that there is no explicit solution for problem \eqref{AGDdeltaprox}. However, following the construction of the past section, we can approximate the solution $X^\delta(u,v)$ by solving the fixed-point equation 
$$s = \sum_{ij} \max\left\{\exp\left(\frac{-C_{ij}+u_i+v_j}{\gamma}\right),s\frac{{\delta}}{n^2} \right\}.$$

\subsection{Nonlinear penalization}\label{sec:nlp}
Here we observe that an advantage of the saddle-point formulation and the methods studied in this paper is that they can be recast without effort to address nonlinear variants of Optimal Transport and Wasserstein Barycenter problems. Nonlinear extensions cover interesting problems such as unbalanced OT and WB problems,  which could not be tackled easily with exact linear optimization-based methods.

We start by defining a simple case of nonlinear Optimal Transport. Let $\psi: \mathbb{R}^n \rightarrow \mathbb{R}_+$ a convex function and consider the following penalized OT problem and its saddle-point formulation:
\begin{equation}
\min_{X\in \Delta_\mu} \la C,X\ra + \psi(\nu-X^\top\onen) = \min_{X\in \Delta_\mu}\max_{v} \la C - \onen \otimes v,X\ra +\la v, \nu \ra-\psi^*(v).
\label{penOT}
\end{equation}
The problem above can be solved with the Hybrid Primal Dual method as soon as the proximal problem 
\begin{equation}
\textup{arg}\min_{v} \langle v,X^\top\onen - \nu\rangle + \psi^*(v)+\frac{1}{2\tau}\|v-\bar v\|^2_2,
\label{prox-pen}
\end{equation}
can be solved efficiently. Note that the proximal problem associated to $X$ does not depend on the penalization $\psi$, therefore the solution does not change (see Sec.~\ref{sec:non-regOT}).

Formulation \eqref{penOT} can be used to deal with unbalanced Optimal Transport, i.e. an Optimal Transport problem where the two marginals have different total mass (\textit{cf}~for
instance \cite[Sec.~10.2]{PeyreCuturibook}).
In this setting, the classical OT problem \eqref{OT} does not admit a feasible solution,
and one looks for a transportation matrix that satisfies the marginal constraints approximately.
(In addition, we may assume that $\mu \cdot \onen = 1$ so that $X$
still belongs to the simplex --- if not, it belongs to a simplex multiple of the unit simplex
and our analysis is easy to adapt).

Several penalization functions can be used to approximate the constraint $\nu - X^\top \onen= 0$,
see again~\cite{PeyreCuturibook}.
For instance, defining $\psi= \frac{1}{2\eta}\|\cdot\|^2_2$ leads to quadratic penalization balanced with the parameter $\eta>0$. Remark that when $\eta \rightarrow 0$, problem \eqref{penOT} becomes equivalent to the classical formulation \eqref{OT}. Moreover,  for a quadratic penalization,
\eqref{prox-pen} has an explicit solution  given by $v = \frac{\eta}{\eta+\tau}\left(\bar v + \tau(\nu-X^\top\onen)\right)$.
One may also
consider a total variation penalization, that is $\psi = \alpha\|\cdot\|_1$.
It leads to the explicit update $v = \Pi_{[-\alpha,\alpha]}\left(\bar v + \tau(\nu-X^\top\onen) \right)$. Note that if $\alpha = \lambda = \|C\|_\infty/2$, this is already the update of equation \eqref{solutionsOT}. 

The same ideas can also be extended to the setting of unbalanced Wasserstein Barycenter. In that case, we may consider the saddle-point formulation 
\begin{equation}
\min_{\nu\ge 0,\;  X^l\in {\Delta_\nu}} \max_{v} \sum_{l=1}^m w_l\left(\langle C_l-\onen \otimes v_l,X_l\rangle +\langle v_l,\nu\rangle - \psi_l^*(v_l)\right),
\end{equation}
where $\psi_l$ are the penalization functions. Then, we recover explicit updates  for all $l=1,...,m$
\begin{equation*}
X^{k+1}_l = \mathcal{N}_\mu\left(X^k_l\exp \left(-\sigma_k(C_l - \onen\otimes \bar v^k_l\right)\right), \qquad \nu^{k+1}_l = \nu^k\exp\left( -\sigma_k\sum_{l=1}^m w_l\bar v_l\right),
\end{equation*}
and $v^{k+1}_l = \frac{\eta}{\eta+\tau_k}\left(v^k_l + \tau_k(\nu^{k+1}-(X^{k+1}_l)^\top\onen)\right)$ in the case of quadratic penalization, or  $v^{k+1}_l = \Pi_{[-\eta,\eta]}\left( v^k_l + \tau_k(\nu^{k+1}-(X^{k+1}_l)^\top\onen) \right)$ for total variation penalization. 

An example of this using quadratic penalization is computed in Figure~\ref{fig:plot6} (bottom right). Observe however that a specific complexity analysis should be performed here, depending on the nonlinearity (as for instance we cannot use straight out of the box a result such as lemma~\ref{thm:halfC} in this particular setting), which is out of the scope of the present study.

\section{Numerical Experiments}
\label{sec:experiments}

In this section, we present numerical results over different instances of Optimal Transport and Wasserstein Barycenter problems for the algorithms described in this paper and for methods in the literature. We implement the methods using \verb'Python-torch' and run the algorithms on a MacBook Air M1 Octa Core with 8GB RAM. We thank authors in \cite{guminov2021combination} and \cite{dvinskikh2021improved} for sharing their implementations in \verb'Python-numpy', which we replicate in \verb'Python-torch' in order to have a fair comparison. 

\subsection{Instances}
For Optimal Transport problems we consider the following instances. 
\paragraph{MNIST instance:} we randomly choose two images from the MNIST library, normalize and use it as the marginals $\mu$ and $\nu$. Note that images in the MNIST library are handwritten numbers of size $28\times 28$ pixels, hence $n = 784$ in this case. For the cost matrix we stick to the Euclidean distance setting.

\paragraph{Gaussian instance:} we measure the OT distance between distributions $\mu \sim \mathcal{N}(3,1)+\mathcal{N}(7,1)$, and $\nu \sim \mathcal{N}(5,1)$. To generate a discretization of this problem, we consider a partition of the interval $[0,10]$ into $n$ equidistant points, then use the Euclidean distance as a metric. 

\paragraph{Random instance:} given $n>0$, we draw the components of the vectors $\mu,\nu \in \mathbb{R}^n$, and the matrix $C\in\mathbb{R}^n\times \mathbb{R}^n$ independently from an uniform distribution in $[0,1]$. 

\begin{figure}[ht]
\centering
\includegraphics[scale=.27]{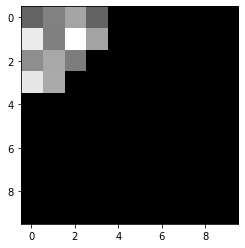}
\includegraphics[scale=.27]{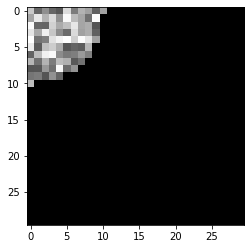}
\includegraphics[scale=.27]{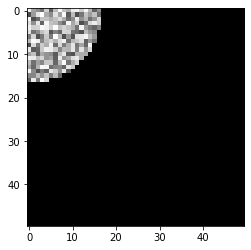}
\includegraphics[scale=.27]{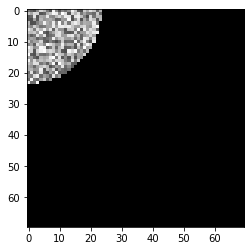}
\includegraphics[scale=.27]{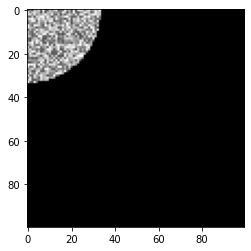}
\includegraphics[scale=.27]{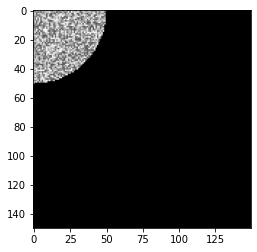}
\caption{Corner images for  $n\in \{10,30,50,70,100,150\}$}
\label{corner_to_dense}
\end{figure}

\paragraph{Corner to dense:} in these instances, the first marginal correspond to an image of size $n_{\textup{pix}}\times n_{\textup{pix}}$ where the source is concentrated in the top-left corner (see Figure~\ref{corner_to_dense}). As a second marginal, we consider an image of the same size with uniform distribution of the mass. We use these instances to measure the performance of exact solvers when transporting the mass from a sparse into a dense marginal. We remark that for those instances the marginals are unbalanced, therefore we normalize it so the total mass sum up to one.

\subsection{Exact solutions}
In order to compute exact solutions for the instances, we use the Python Optimal Transport (POT) library \cite{flamary2021pot}. In particular, for optimal transport we use the  \verb'ot.emd' function that is based on the algorithm in \cite{bonneel2011displacement}, while for Wasserstein Barycenter problems we use the  \verb'ot.lp.barycenter' function which is based on Network Simplex. Table~\ref{exact_solver} shows the performance of the exact solver and our approximated method for different instances and sizes of Optimal Transport problems.

\begin{table}[ht]
\centering
\begin{tabular}{|ccc|ccc|ccc|}
\hline
\multicolumn{9}{|c|}{Exact solver} \\ \hline
\multicolumn{3}{|c|}{Corner to dense} & \multicolumn{3}{|c|}{Gaussian} & \multicolumn{3}{|c|}{Random} \\ 
$n_{\textup{pix}}$ & time $(s)$ & $\ell^1$ error & $n$ & time $(s)$ & $\ell^1$ error & $n$ & time $(s)$ & $\ell^1$ error \\ \hline 
10 & 0.0038 & 1.10e-07 & 100 & 0.0011 & 3.27e-16 & 100 & 0.0015 & 4.18e-16 \\ 
30 & 0.5275 & 1.15e-07 & 500 & 0.0238 & 6.90e-16 & 500 & 0.0227 & 6.97e-16 \\
50 & 10.5236 & 2.64e-07 & 1000 & 0.1495 & 1.02e-15 & 1000 & 0.1012 & 7.83e-16 \\
70 & 92.6429 & 8.84e-08 & 5000 & 3.7442 & 1.97e-15 & 5000 & 4.8426 & 1.12e-15 \\
100 & - & - & 10000 & 148.89 & 2.39e-15 & 10000 & 87.1855 & 1.60e-15 \\ 
120 & - & - & 15000 & - & - & 20000 & - & - \\ \hline \hline
\multicolumn{9}{|c|}{Regularized Hybrid Primal Dual with linesearch} \\ \hline
\multicolumn{3}{|c|}{Corner to dense} & \multicolumn{3}{|c|}{Gaussian} & \multicolumn{3}{|c|}{Random} \\ 
$n_{\textup{pix}}$ & time $(s)$ & precision & $n$ & time $(s)$ & precision & $n$ & time $(s)$ & precision \\ \hline 
10 & 0.0017 & 0.0089 & 100 & 0.0554 & 0.0095 & 100 & 0.1600 & 0.0098 \\ 
30 & 4.0481 & 0.0064 & 500 & 0.5636 & 0.0097 & 500 & 0.3982 & 0.0069 \\
50 & 30.0447 & 0.0099 & 1000 & 1.8527 & 0.0091 & 1000 & 0.7227 & 0.0059 \\
70 & 396.647 & 0.0097 & 5000 & 119.69 & 0.0098 & 5000 & 21.9500 & 0.0046 \\
100 & 558.39 & 0.0215* & 10000 & 139.52 & 0.0175 & 10000 & 274.62 & 0.0038 \\ 
120 & 1437.33 & 0.0230* & 15000 & 502.60 & 0.0201* & 20000 & 494.80 & 0.0021*\\ \hline
\end{tabular}
\caption{Exact solver versus the main algorithm proposed in this paper. The precision measures the distance to the optimal value. For the values with (*), the precision corresponds to the duality gap. }
\label{exact_solver}
\end{table}

We remark that solvers in the POT library are extremely fast to solve medium size instances (see Table~\ref{exact_solver}). However, for large instances ($n>10000$), we did not obtain a solution after 30 minutes of computing time. Note that exact solvers like those in POT are not designed to deal with very large problems (see the POT documentation for extra details). Therefore, approximate methods such as those described in this paper can be of practical interest whenever $n$ is large. 

\subsection{Results for Optimal Transport instances}
\begin{figure}[ht]
\centering
\includegraphics[scale=0.34]{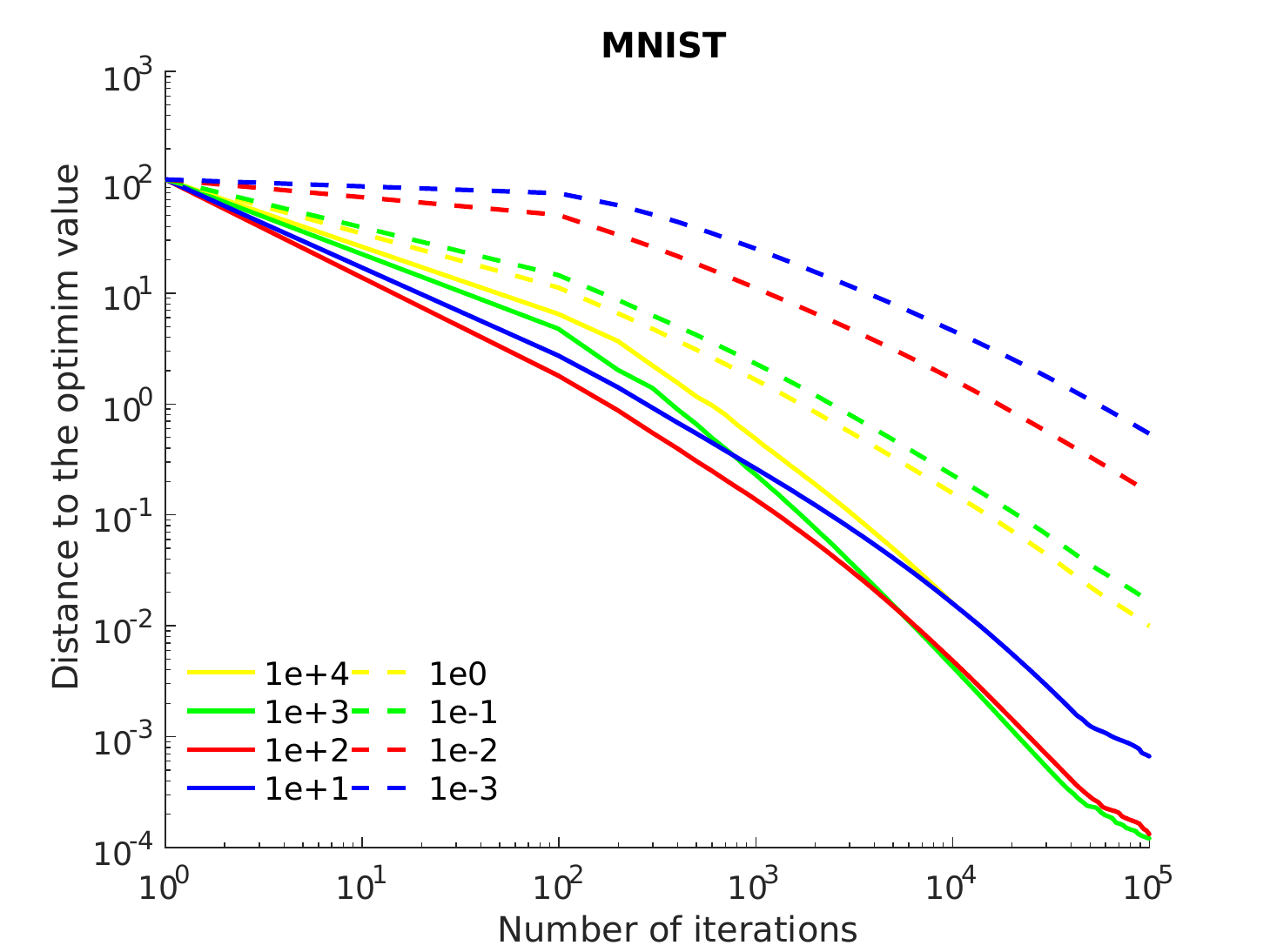} 
\includegraphics[scale=0.34]{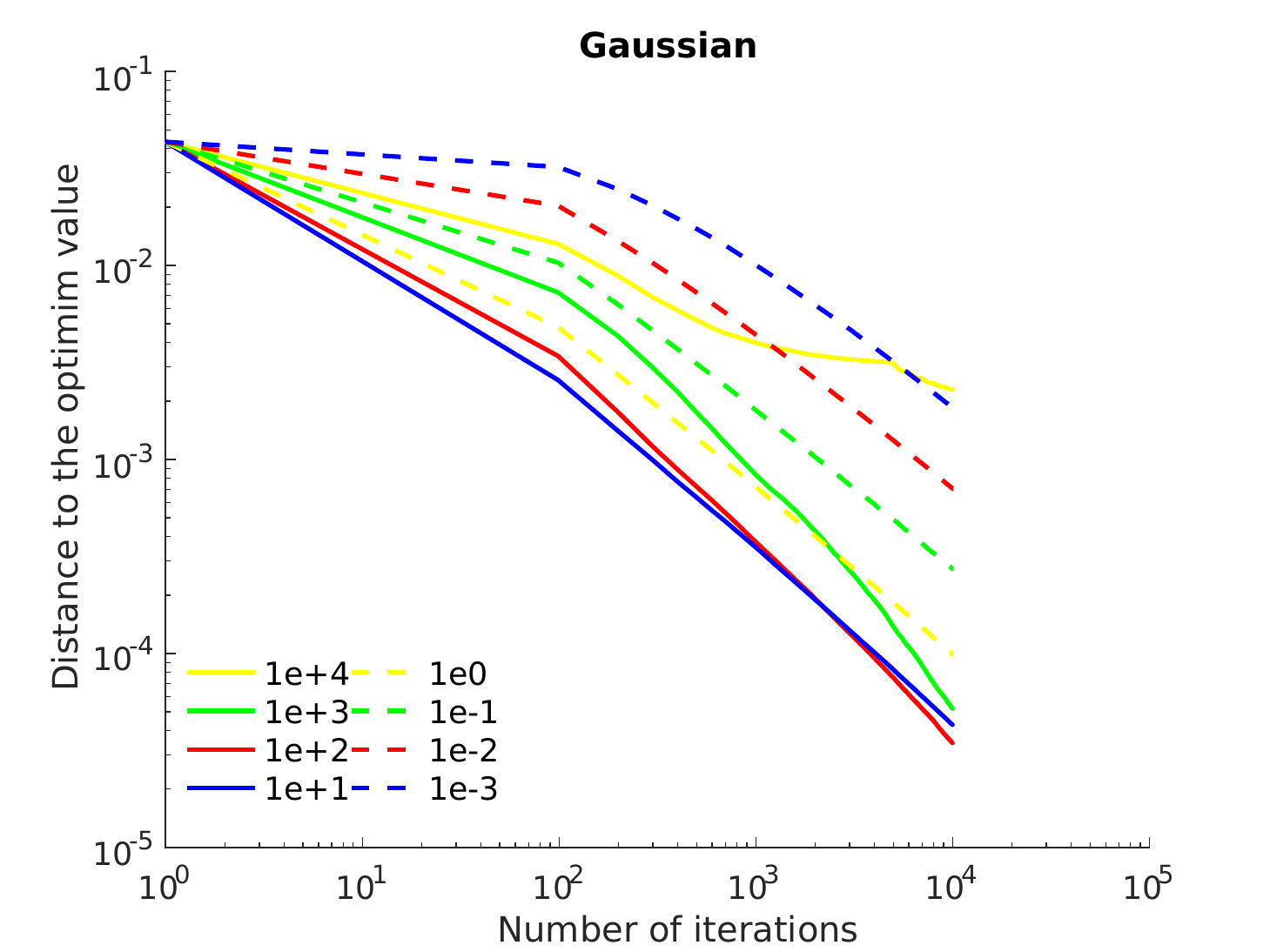} 
\includegraphics[scale=0.34]{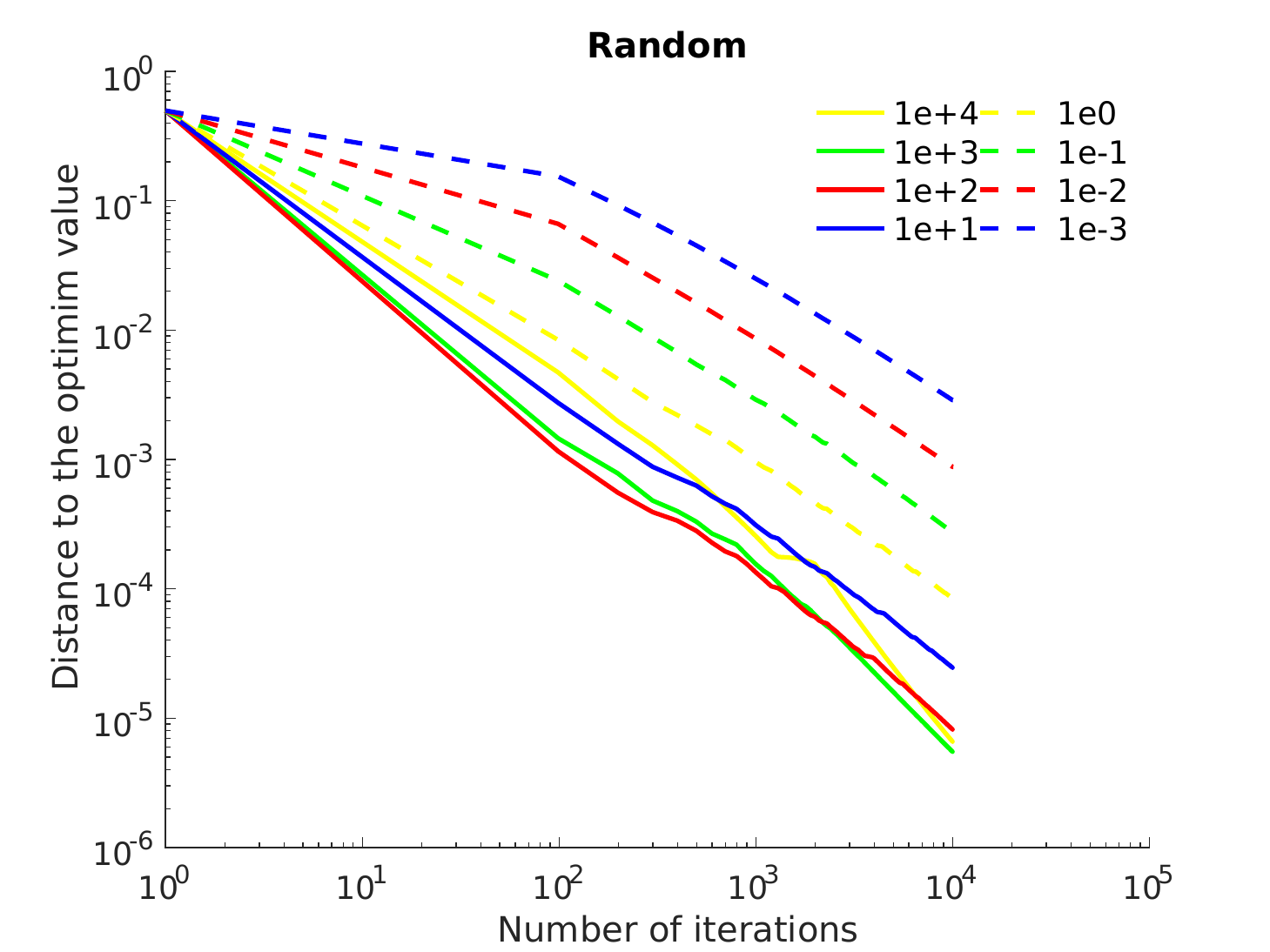} 
\caption{log-log distance to the optimum for {\nothing $\gamma-$HPD\,ls\,fm} with different values $\beta_1 = C\cdot \frac{2\ln n}{n\lambda^2}$, with $C\in\{10^{-3}, 10^{-2}, 10^{-1},1,10, 10^{2},10^{3}\}$. For Gaussian and random instances we set $n=1000$.}
\label{fig:plot0}
\end{figure}

In our first experience, we test the effects of the parameter $\beta_1$ in the performance of the Hybrid Primal Dual method ({HPD}),  with entropy regularization $(\gamma)$, linesearch ({ls}), and the fixed marginals ({fm}) approach discussed in Remark~\ref{fixmarginal}. We refer to this method as {\nothing $\gamma-$HPD\,ls\,fm}. As we discussed in the proof of Proposition~\ref{thm:boundaccel}, there is a compromise between the size of the initial stepsize $\tau_1$ (which decrease with $\beta_1$ increasing), and the bound \eqref{AMP2-OT} (which increase with $\beta_1$ increasing). Figure \eqref{fig:plot0} shows the performance for different choices of $\beta_1 = O(\frac{\ln n}{n\lambda^2})$. We observe that the quality of the solution slightly improves if we amplify the quantity $\frac{\ln n}{n\lambda^2}$, meaning that it is of convenience to start with a smaller initial stepsize $\tau_1$, improving the bound \eqref{AMP2-OT}. Despite this, the optimal value for $\beta_1$ it seems to depend on the characteristic of the problem ($\beta_1 \approx 10^3\frac{\ln n}{n\lambda^2}$ for MNIST and random instances, while $\beta_1 \approx 10^2\frac{\ln n}{n\lambda^2}$ for Gaussians instances). A theoretical study concerning a better choice of $\beta_1$ than provided by the worse case estimates remains an open question.

\begin{figure}[ht]
\centering
\includegraphics[scale=0.34]{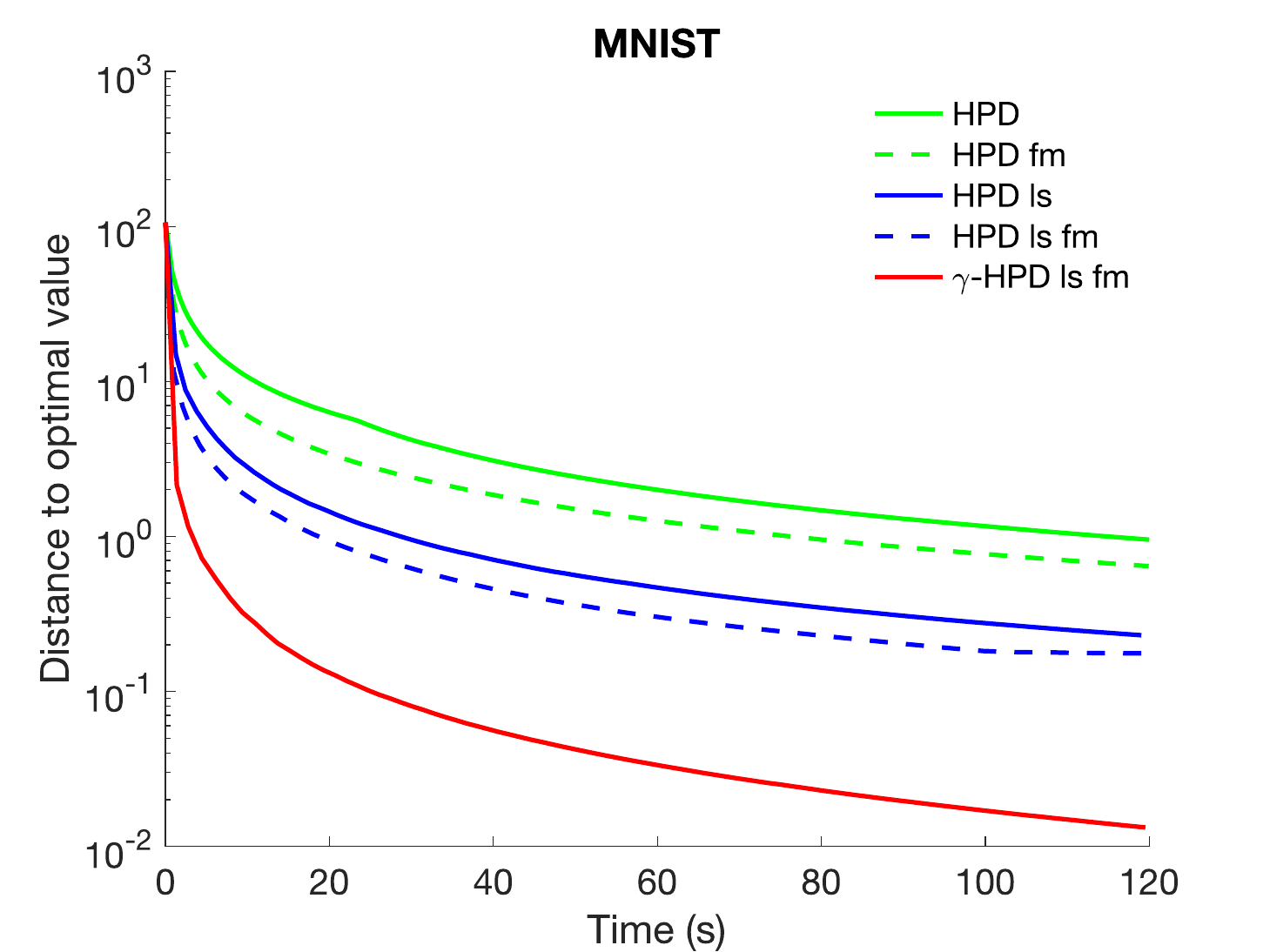} 
\includegraphics[scale=0.34]{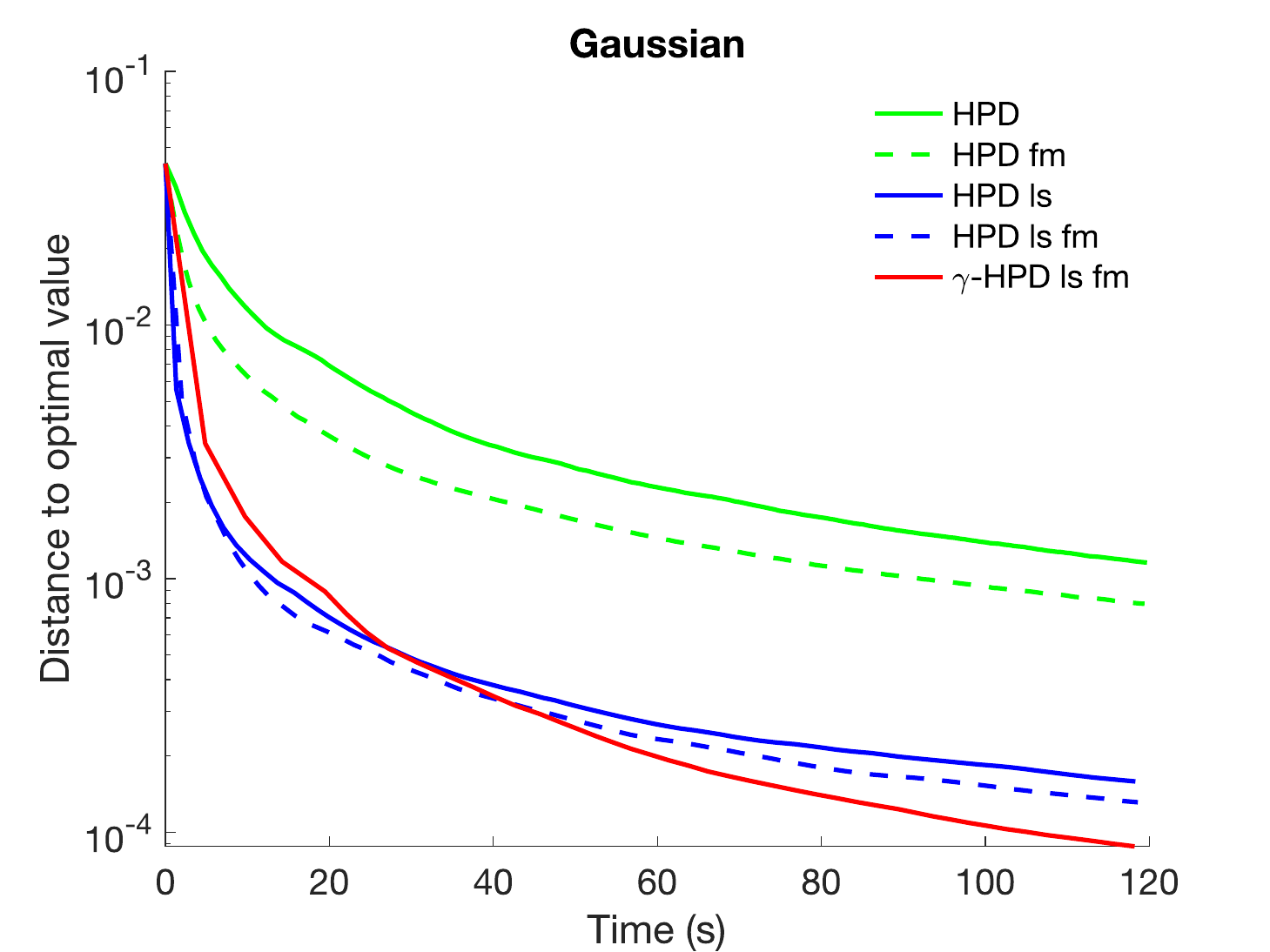} 
\includegraphics[scale=0.34]{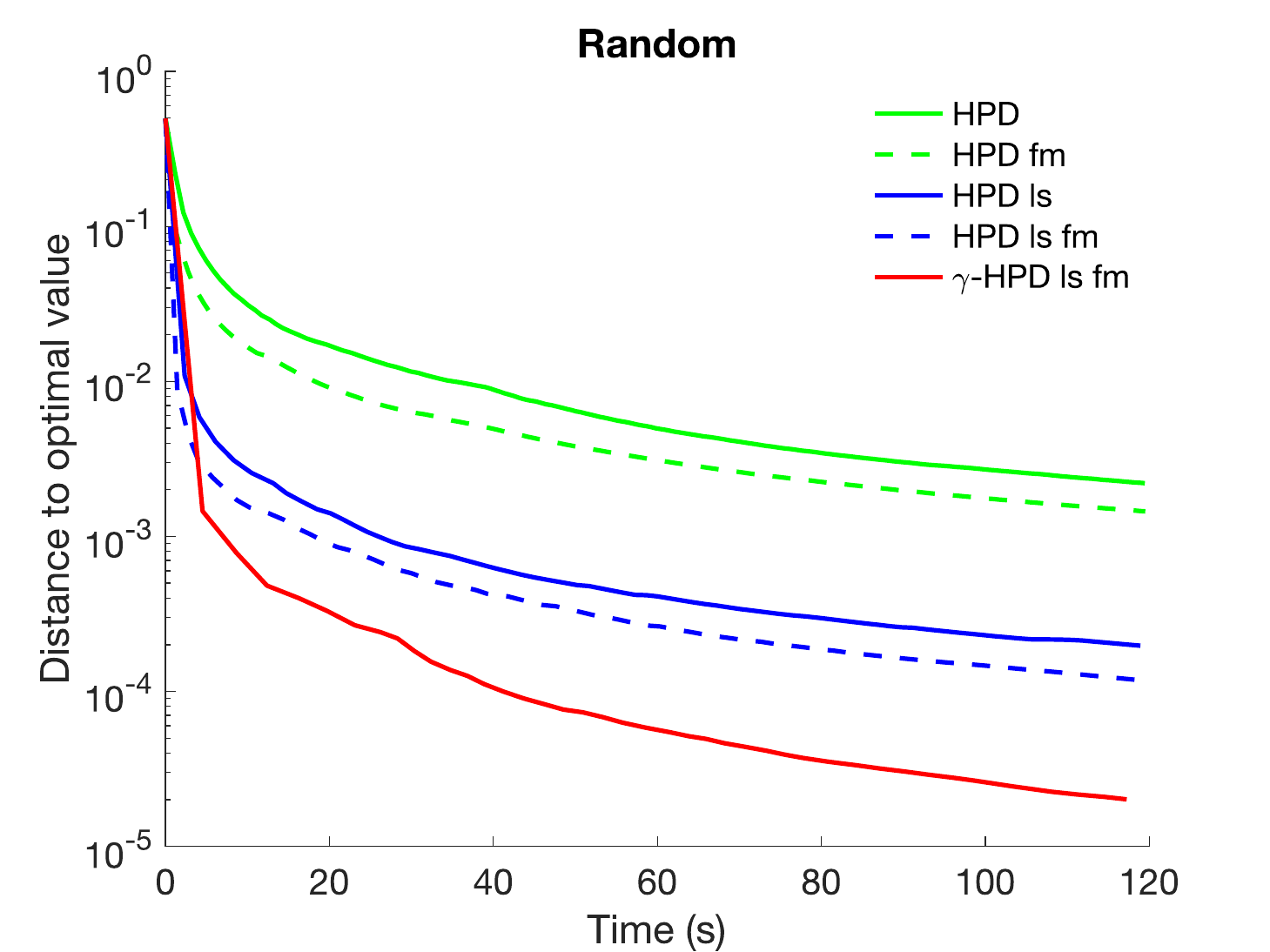} 
\caption{log-log distance to the optimum for different versions of the HPD method, including linesearch ({\nothing ls}), fixed-marginals {(\nothing fm)}, and entropic regularization ($\gamma$).}
\label{fig:plot0-1}
\end{figure}

Next, we compare the performance of several versions of the HPD method including: non-regularized {\nothing (HPD)}, entropy regularized ({\nothing $\gamma-$HPD}), and its versions with linesearch ({\nothing ls}) and fixed marginals ({\nothing fm}). Figure~\ref{fig:plot0-1} shows that, for the three instances, the best performance is attained when we include linesearch and fixed marginals. Moreover, we observe that using fixed marginals also improves the performance of non-regularized versions of the method.

\begin{figure}[ht]
\begin{center}
  \includegraphics[scale=0.34]{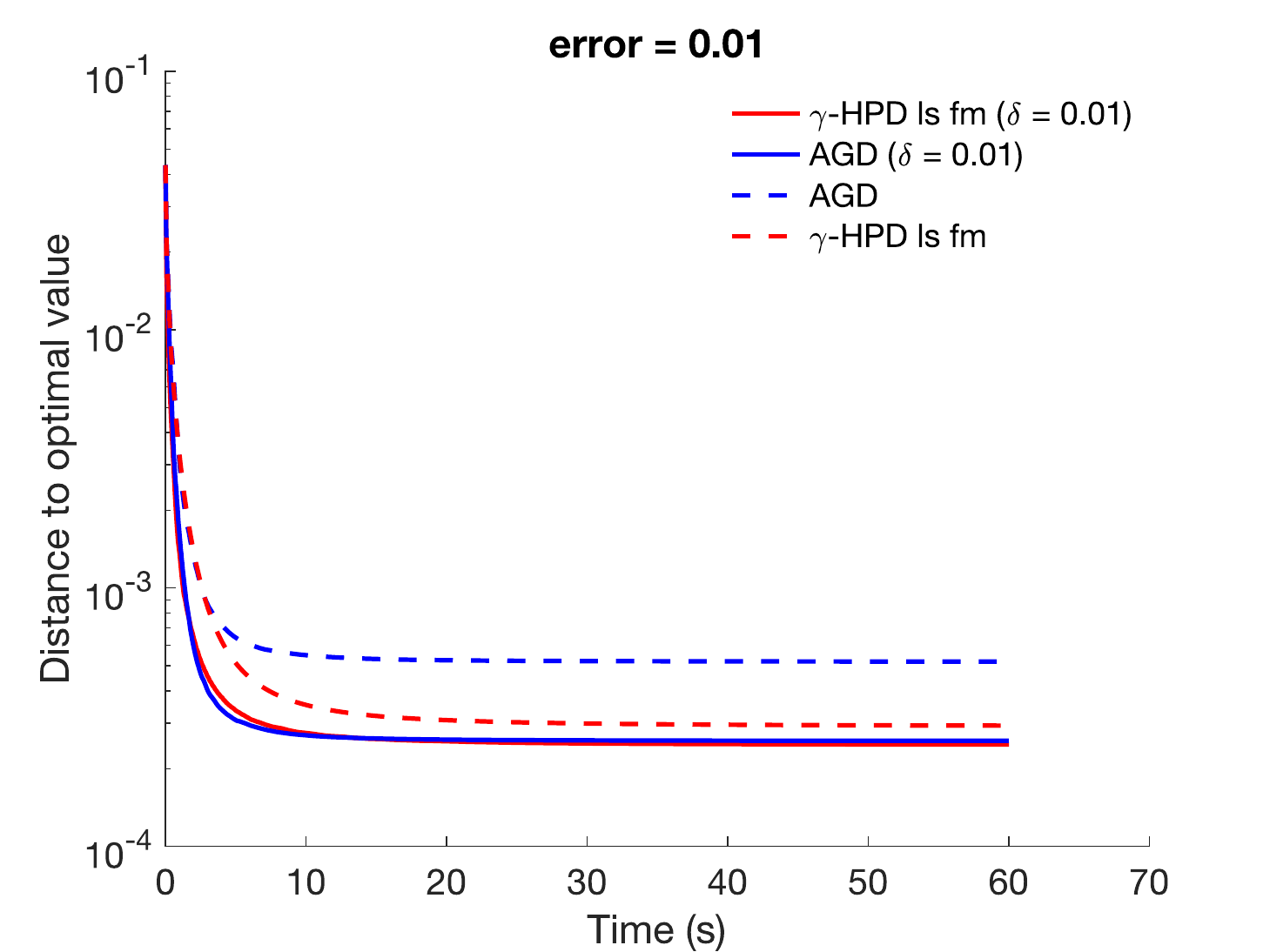}
  \includegraphics[scale=0.34]{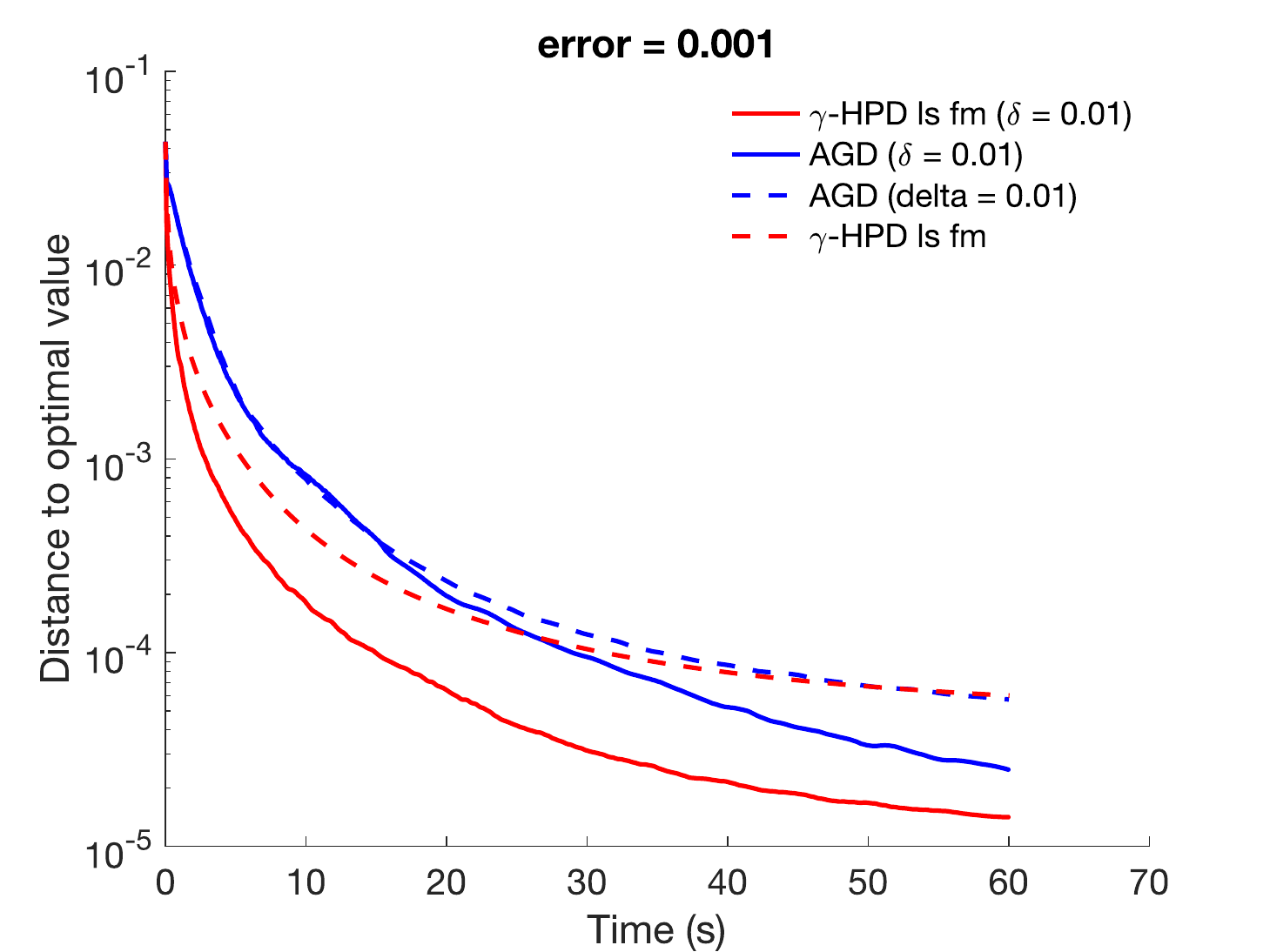}
  \includegraphics[scale=0.34]{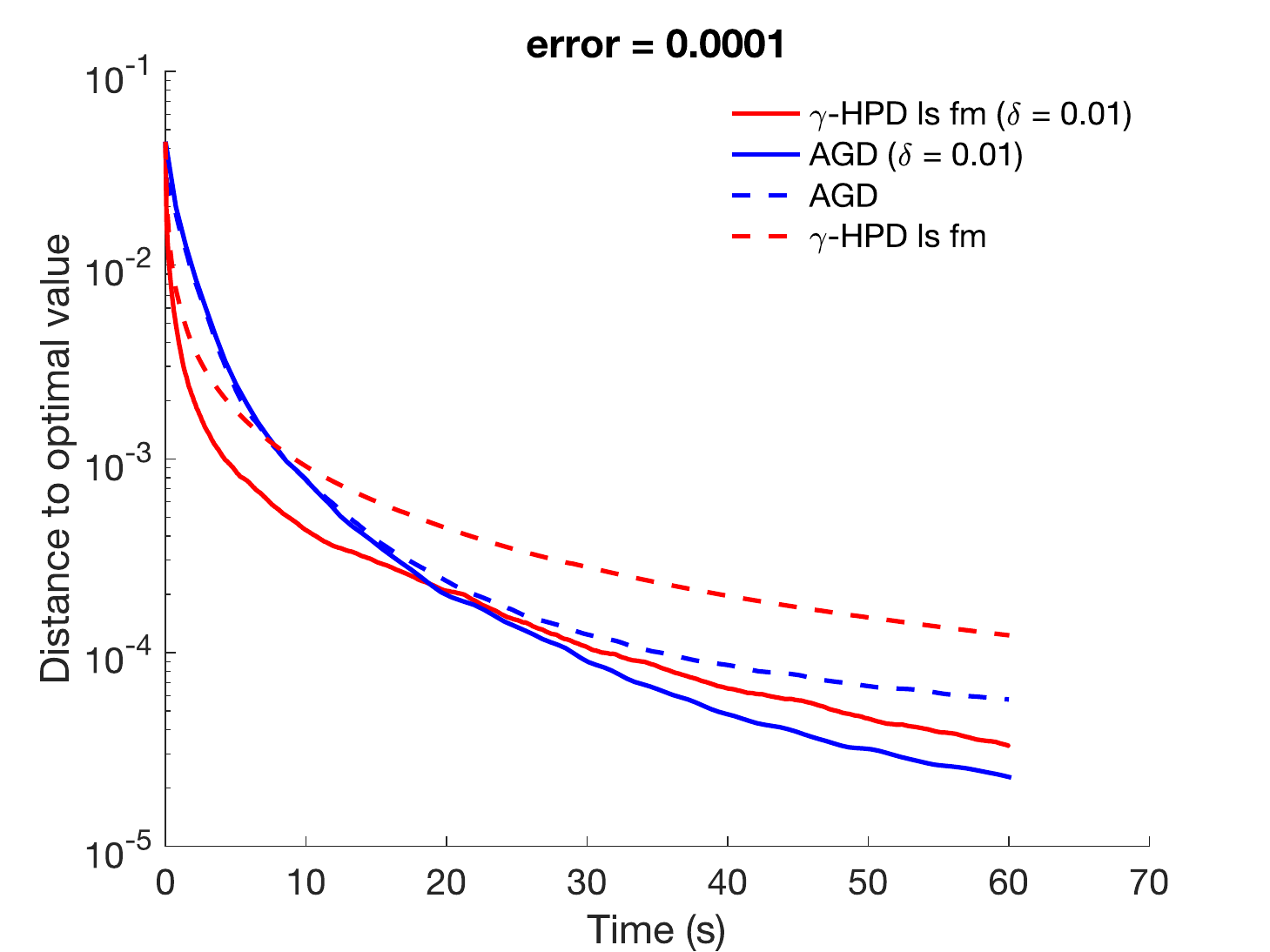}
\caption{Distance to the optimum value for different prefixed tolerances.}
\label{fig:plot2}
\end{center}
\end{figure}

\begin{figure}[ht]
\begin{center}
\includegraphics[scale=1]{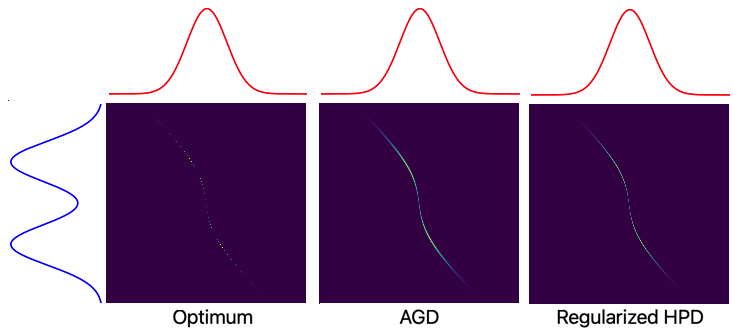}
\end{center}
\caption{Transportation plan between Gaussian distributions reported by {\nothing $\gamma-$HPD\,ls} and {\nothing AGD} with $delta$-regularization ($\delta = 0.01$).}
\label{fig:plot_delta-sol}
\end{figure}

In the third experiment, we test the impact of the prefixed parameter $\varepsilon>0$ in the performance of our algorithm {\nothing $\gamma$-HPD\,ls\,fm}, as well as the Accelerated Gradient Descent ({\nothing AGD}), both with and without the $\delta-$regularization of Section~\ref{sec:sek} (see Figure~\ref{fig:plot2}). We observe that {\nothing $\gamma$-HPD\,ls\,fm} is very stable when the tolerance approaches zero. We remark that the sensibility to $\varepsilon$ highly depends on the implementation strategy to avoid numerical issues. An extensive comparison of different implementation strategies escapes the scope of this paper and we leave it as possible future work. Finally, the $\delta-$regularized methods are also stable and perform slightly better than non-regularized methods. Moreover, in Figure~\ref{fig:plot_delta-sol}, we show the structure of the optimal transport reported by {\nothing $\gamma-$HPD\,ls} and {\nothing AGD} when using the scaled entropy regularization. We remark that these solutions are sparse and have a support similar to the exact optimal solution, contrarily to solutions obtained with entropy regularized methods, which are never sparse.

\subsection{Results for Wasserstein Barycenter problems}

\begin{figure}[ht]
\centering
\includegraphics[scale=0.25]{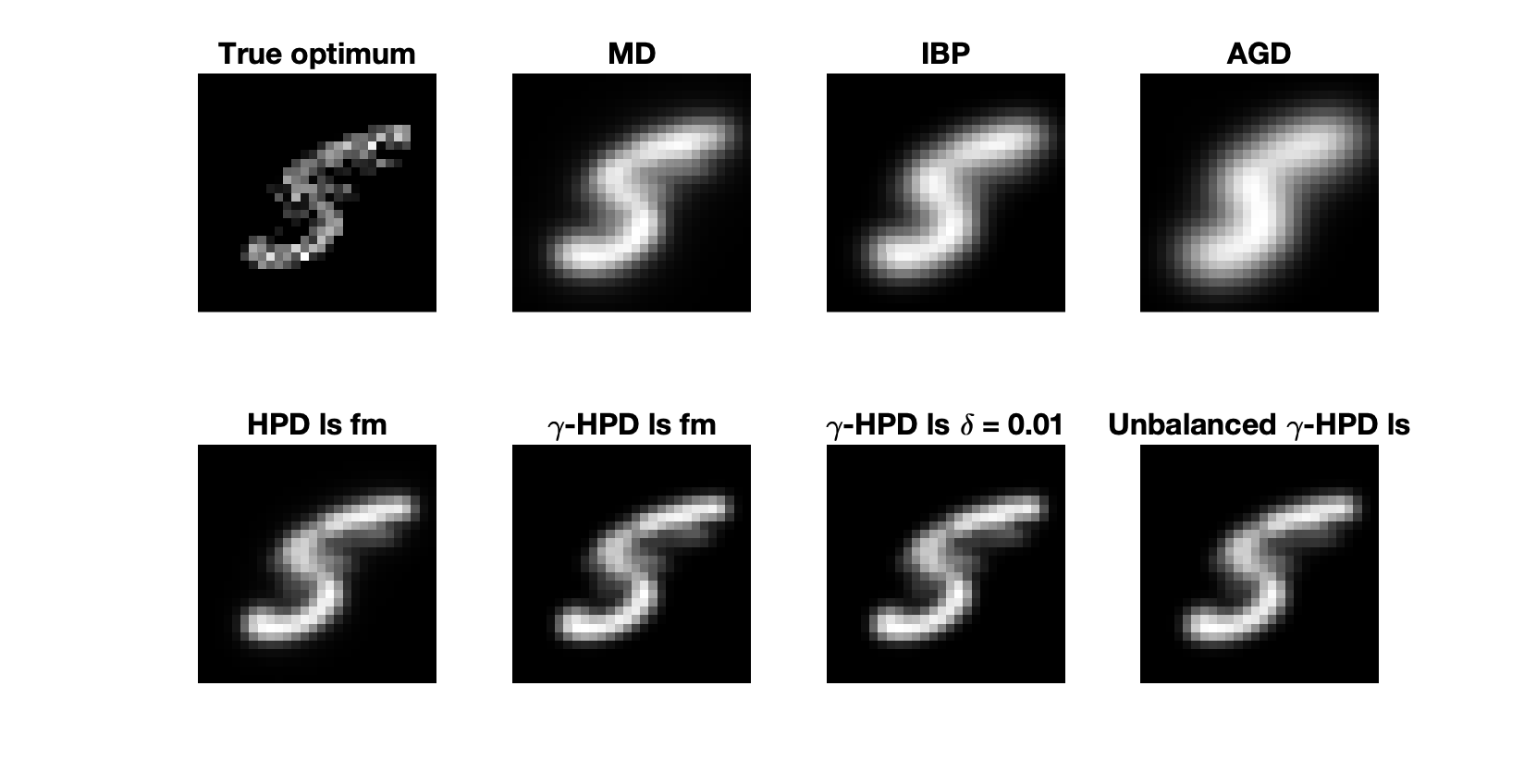}
\caption{Barycenters after 10000 iterations of each method.}
\label{fig:plot6}
\end{figure}

\begin{figure}[ht]
\centering
\includegraphics[scale=0.7]{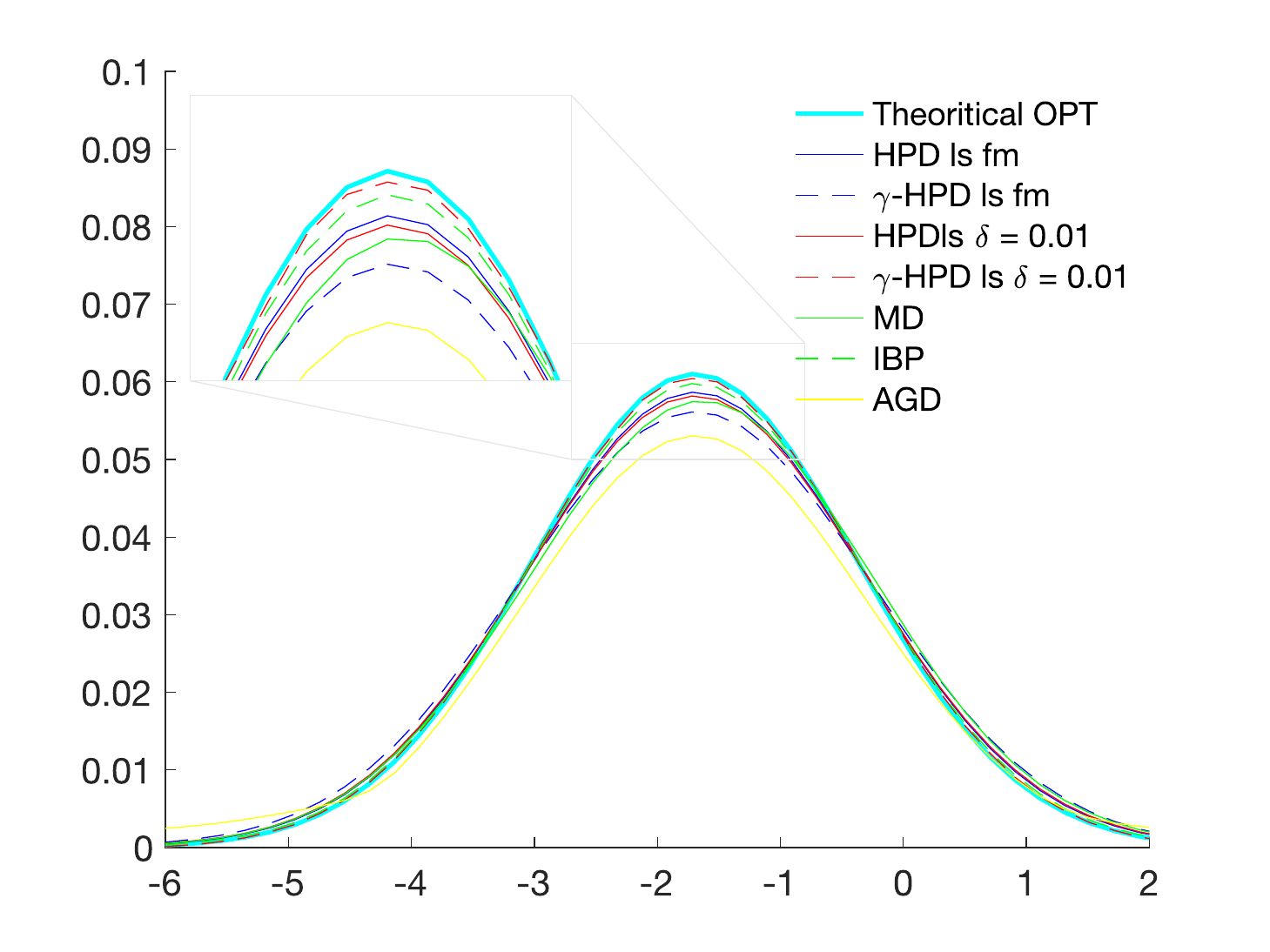} 
\caption{Wasserstein Barycenters computed after 10000 iterations of each method.}
\label{fig:plot7}
\end{figure}

For Wasserstein Barycenter problems we first use a collection of $m=5$ images of the number five in the MNIST library. Figure \ref{fig:plot6} shows the barycenter computed for each method. We also plot the true optimum at the top-left for reference. Then, we compare the performance of the methods when solving a Wasserstein Barycenter problem between 1-dimensional Gaussian measures. For this, we generate $m=10$ different Gaussian measures with random parameters and $n=100$ bins in the interval $[-10,10]$. Figure \eqref{fig:plot7} shows the barycenters obtained with the algorithms of literature and the ones presented in this paper. We also plot the theoretical barycenter computed according to Theorem~2.2 in \cite{ruschendorf2002n}.

The results for WB show that our proposed algorithms are very competitive with the approaches in the literature. In this case, the solutions with $\delta-$regularization are slightly sharper yet not as sparse as for OT. An explanation is that the support of the barycenter results from a mixture of the different supports of multiple transportation matrices, and should therefore have a larger support. Finally, we mention that the bottom-right image in Figure \ref{fig:plot6}, which is also sharp, corresponds to the regularized HPD applied directly to the unbalanced MNIST marginals and using a quadratic penalization as described in Section~\ref{sec:nlp}

\section{Conclusions}

We have shown how to apply and extend first-order primal-dual methods for saddle-point problems in a non-Euclidean setting to approximate the
solutions of Optimal Transport and Wasserstein Barycenter problems. Furthermore, the proposed methods provide a new perspective between the Mirror-Prox and the AGD as they may or may not rely on entropic regularization and use linesearch to improve the numerical performance. We believe that this new perspective is of great interest to understand the dynamics behind the accelerated OT methods and may help tackle open questions, such as the existence of a competitive $O(n^2/\varepsilon)$ algorithm. Afterward, we introduced a modified entropy kernel for the Bregman mirror descent steps that improve the solution's quality at a low computational cost. This idea can also be accomplished with AGD-type methods to obtain sparse transportation plans that are of interest in applications and cannot be obtained with the traditional approach. Numerical experiments show that the presented algorithms are very stable and competitive with the state-of-art-methods in time and precision.

\section*{Acknowledgements}
The work of Juan Pablo Contreras was supported by a doctoral scholarship from ANID-PFCHA/Doctorado Nacional/2019-21190161. The authors warmly thank the reviewers of this manuscript for their comments and their very helpful suggestions.



\appendix

\section{Proofs of the main results}


\subsection{Detailed proof of Theorem \ref{thm:MP}}\label{sec:proofMP}

We mostly focus on the difference with the proofs in~\cite{malitsky2018first}.
In particular, the termination criterion we consider in Alg.~\ref{Alg:PDls}
in Step~\ref{MPrule} terminates earlier than the one in~\cite{malitsky2018first}
and it can be shown that it induces a globally bounded (actually,
depending on $\tau_0,\beta_0,\rho,L$ in a logarithmic way)
extra multiplicative factor in the complexity.

We start by recalling the definition of a relatively strongly convex function.

\begin{definition} \cite{Bauschkeetal-nolips,lu2018relatively}
Given $\gamma>0$, the function $f$ is $\gamma-$strongly convex relative to $\xi$ if for all $x,y \in \textup{int}\,\textup{dom}\,\xi$ it holds
$$f(y)\ge f(x)+\langle \nabla f(x),y-x\rangle+\gamma D_\xi(y,x).$$
In particular $f$ is $\gamma$-strongly convex relative to $\xi$
if and only if  $f-\gamma\xi$ is convex over $\textup{int}\,\textup{dom}\,\xi$.
\end{definition}

Let us recall the following results, the first being a particular case of the second.
  \begin{lemma}[Proximal optimality condition]
  Let $f:\mathcal{X}\rightarrow \mathbb{R}\cup\{+\infty\}$ a closed, convex function.
      Then $\hat x$ is a solution of the proximal problem $\min_{x\in \mathcal{X}} f(x)+\frac{1}{2}\|x-\bar x\|^2_2$, if and only if
$$ \langle \hat x-\bar x,\hat x- x \rangle + f(x) \ge f(\hat x), \qquad \forall x\in \mathcal{X}.$$
\label{proxopt}
\end{lemma}

\begin{lemma}[Chen and Teboulle 1993]\label{lem:3p}
Let $f:\mathcal{X}\rightarrow \mathbb{R}\cup\{+\infty\}$ be a closed proper convex function and $x^+$ an optimal solution to the nonlinear prox problem $\displaystyle{\min_{x\in \mathcal{X}}\left\{f(x)+\tau^{-1}D_\xi(x,\bar x)\right\}}$. Then, for any $x\in {\textup{dom}\,\xi}$ it holds
$$\tau(f(x^+)-f(x))\le D_\xi(x,\bar x)-D_\xi(x,x^+)-D_\xi(x^+,\bar x).$$ 
\end{lemma}

\paragraph{Proof of Theorem \ref{thm:MP}.} Let $x\in \textup{dom}\,g$ and $y\in \textup{dom}\, \xi_{\mathcal{Y}}$, where we recall that $\xi_\mathcal{Y}$ is a Legendre function (see \cite{Rockafellar}). Remark that in the following analysis $(x,y)$ is not necessarily a saddle-point. We mostly follow~\cite{malitsky2018first}.

Because of the convexity of $g$ and the proximal optimality condition in Lemma~\ref{proxopt}, we have
\begin{equation}
\langle x^{k}-x^{k-1}+\tau_{k-1}K^*y^{k},x-x^{k}\rangle \ge \tau_{k-1}(g(x^{k})-g(x)).
\label{AThm1:1}
\end{equation}
Taking $x=x^{k+1}$ in \eqref{AThm1:1}  multiplied by $\theta_k=\frac{\tau_k}{\tau_{k-1}}$, plus \eqref{AThm1:1} with $x=x^{k-1}$ multiplied by $\theta_k^2$, we obtain
\begin{equation}
\langle \bar x^k-x^k+\tau_{k}K^*y^k, x^{k+1}-\bar x^k\rangle \ge \tau_{k}((1+\theta_k)g(x^k)-g(x^{k+1})-\theta_kg(x^{k-1})).
\label{AThm1:2}
\end{equation}
Using the 3-point inequality of Lemma~\ref{lem:3p} with $h^*$ and $g$ gives, respectively 
\begin{equation}
D_{\mathcal{Y}}(y,y^k)-  D_{\mathcal{Y}}(y^{k+1},y^k) -(1+\gamma)D_{\mathcal{Y}}(y,y^{k+1}) \ge 
 \sigma_k(h^*(y^{k+1})-h^*(y)+\langle K\bar x^k,y-y^{k+1}\rangle). 
 \label{Thm1:3}
\end{equation}
\begin{equation}
\|x-x^k\|^2_2 -  \|x-x^{k+1}\|^2_2 -\|x^{k+1}-x^k\|^2_2 \ge  
 2\tau_{k}(g(x^{k+1})-g(x)+\langle Kx^{k+1}-Kx,y^{k+1}\rangle). 
 \label{Thm1:4}
\end{equation}
Above we assume $h^*$ is $\gamma-$strongly convex relative to $\xi_{\mathcal{Y}}$. If $h^*$ is just convex, then the inequalities holds with $\gamma = 0$.

Summing \eqref{AThm1:2}, \eqref{Thm1:3} and \eqref{Thm1:4} and in view of the stopping criteria for the linesearch (step~\ref{MPrule} in Algorithm~\ref{Alg:PDls}),
reproducing the same computation as~\cite{malitsky2018first} (but
here $(x,y)$ is not a saddle-point and we cannot rely on the
non-negativity of expressions such as $\mathcal{L}(x^k,y)$):
\begin{align}
(1+\theta_k)&\mathcal{L}(x^k,y)-\theta_k\mathcal{L}(x^{k-1},y)-\mathcal{L}(x,y^{k+1})\le \nonumber \\
& \frac{1}{2\tau_k}\|x^k-x\|^2_2-\frac{1}{2\tau_k}\|x-x^{k+1}\|^2_2+ \frac{1}{\sigma_k}D_{\mathcal{Y}}(y,y^k)-\frac{1+\gamma \sigma_k}{\sigma_k}D_{\mathcal{Y}}(y,y^{k+1}).
\label{Thm:8}
\end{align}
Let $r_k = (1+\theta_k)\mathcal{L}(x^k,y)-\theta_k\mathcal{L}(x^{k-1},y)-\mathcal{L}(x,y^{k+1})$. Equivalently equation \eqref{Thm:8} can be written as
\begin{multline}
r_k \le  \frac{1}{2\tau_k}\|x^k-x\|^2_2  -\frac{\tau_{k+1}}{\tau_k}\frac{1}{2\tau_{k+1}}\|x-x^{k+1}\|^2_2  \\
+ \frac{1}{\sigma_k}D_{\mathcal{Y}}(y,y^k)-\frac{\sigma_{k+1}}{\sigma_k}(1+\gamma \sigma_k)\frac{1}{\sigma_{k+1}}D_{\mathcal{Y}}(y,y^{k+1}).
\end{multline}
By construction of the algorithm, $\frac{\tau_{k+1}}{\tau_k} = \frac{\sigma_{k+1}}{\sigma_k}(1+\gamma\sigma_k)$. Let $R_k := \frac{1}{2\tau_k}\|x-x^k\|^2_2+\frac{1}{\sigma_k}D_\mathcal{Y}(y,y^k)$ then $\tau_kr_k \le \tau_k R_k - \tau_{k+1}R_{k+1}$. Summing the inequalities for $k=1,...,N$ gives 
\begin{equation}
\sum_{k=1}^N \tau_kr_k \le \tau_1R_1 - \tau_NR_N \le \tau_1R_1 =\frac{1}{2}\|x-x^1\|^2_2+\frac{1}{\beta_1}D_\mathcal{Y}(y,y^1)
\end{equation}
On the other hand, we have
\begin{multline*}
  \sum_{k=1}^N \tau_k(1+\theta_k)\mathcal{L}(x^k,y)-\tau_k\theta_k\mathcal{L}(x^{k-1},y)  = \sum_{k=2}^N (\tau_{k-1}(1+\theta_{k-1})-\tau_k\theta_k)\mathcal{L}(x^{k-1},y)
  \\+\tau_N(1+\theta_N)\mathcal{L}(x^N,y)-\tau_1\theta_1\mathcal{L}(x^0,y).
\end{multline*}
Note that $\tau_{k-1}(1+\theta_{k-1})-\tau_k\theta_k\ge 0$ and, recalling $T_N=\sum_{k=1}^N\tau_k$:
\[
  \sum_{k=2}^N \left( \tau_{k-1}(1+\theta_{k-1})-\tau_k\theta_k\right) + \tau_N(1+\theta_N) = \tau_1\theta_1+T_N.
\]
Hence, by convexity of $\mathcal{L}(\cdot,y)$, one has:
\begin{multline*}
  \sum_{k=2}^N (\tau_{k-1}(1+\theta_{k-1})-\tau_k\theta_k)\mathcal{L}(x^{k-1},y)+\tau_N(1+\theta_N)\mathcal{L}(x^N,y)
\\  \ge
(\tau_1\theta_1+T_N)  \mathcal{L}\left(
     \frac{\sum_{k=2}^N(\tau_{k-1}(1+\theta_{k-1})-\tau_k\theta_k)x^{k-1} +
    \tau_N(1+\theta_N)x^N}{\tau_1\theta_1+T_N}
    ,y\right).
\end{multline*}
Recalling $\bar x^k = (1+\theta_k)x^k-\theta_kx^{k-1}$, one has
\[
 \sum_{k=2}^N(\tau_{k-1}(1+\theta_{k-1})-\tau_k\theta_k)x^{k-1} +
    \tau_N(1+\theta_N)x^N
  = \tau_1\theta_1x^0+\sum_{k=1}^N\tau_k\bar x^k.
\]
Hence we have shown:
\begin{align}
\sum_{k=1}^N \tau_k(1+\theta_k) & \mathcal{L}(x^k,y)-\tau_k\theta_k\mathcal{L}(x^{k-1},y) \ge  \nonumber \\
& (\tau_1\theta_1+T_N)\mathcal{L}\left(\frac{\tau_1\theta_1x^0+\sum_{k=1}^N\tau_k\bar x^k}{\tau_1\theta_1+T_N}, y\right)-\tau_1\theta_1\mathcal{L}(x^0,y),
\label{convexP}
\end{align}
motivating the definition of $\hat x^N = \frac{\tau_1\theta_1x^0+\sum_{k=1}^N\tau_k\bar x^k}{\tau_1\theta_1+T_N}$. Observe that if $\mathcal{L}(\cdot, y)$ is linear in its domain (as is the case for OT), then equation \eqref{convexP} is an equality and the RHS
of~\eqref{convexP}
boils down to $T_N\mathcal{L}\left(\frac{1}{T_N}\sum_{k=1}^N \tau_k\bar x^k,y\right)$, motivating a simpler definition
\begin{equation}\label{eq:alternateX}
  \hat x^N = \frac{1}{T_N}\sum_{k=1}^N \tau_k\bar x^k.
\end{equation}

Similarly, using the convexity of $-\mathcal{L}(x,\cdot)$ we have
\begin{equation}
\sum_{k=1}^N \tau_k\mathcal{L}(x,y^{k+1}) \le T_N\mathcal{L}\left(x,\frac{1}{T_N}\sum_{k=1}^N\tau_ky^{k+1}\right),
\label{convexD}
\end{equation}
motivating the definition $\hat y^N = \frac{1}{T_N}\sum_{k=1}^N\tau_ky^{k+1}$. In brief we have 
\begin{small}
\begin{equation}
T_N(\mathcal{L}(\hat x^N,y)-\mathcal{L}(x,\hat y^N))  \le \frac{1}{2}\|x-x^1\|^2_2+\frac{1}{\beta_1}D_\mathcal{Y}(y,y^1)-\tau_1\theta_1(\mathcal{L}(\hat x^N,y)-\mathcal{L}(x^0,y)).
\end{equation}
\end{small}
Again, when $\mathcal{L}(\cdot,y)$ is linear in its domain, we define $\hat x^N$ by~\eqref{eq:alternateX} and obtain the simpler bound:
\begin{equation}
T_N(\mathcal{L}(\hat x^N,y)-\mathcal{L}(x,\hat y^N))  \le \frac{1}{2}\|x-x^1\|^2_2+\frac{1}{\beta_1}D_\mathcal{Y}(y,y^1).
\label{MPlinear}
\end{equation}

{\bf Complexity estimate: the convex case.}
If $\gamma = 0$, then $\beta_k \equiv \beta$ is constant because of the update rule in step~\ref{step:taubeta}. Moreover, one checks that $\theta_k \le \bar \theta := \frac{1+\sqrt{5}}{2}$ (see~\cite{malitsky2018first}), and $\tau_k \ge \frac{\rho}{\sqrt{\beta}L}$, provided that $\tau_0\ge \frac{\rho}{\sqrt{\beta}L}$, which is the case when we choose $\tau_0 = \frac{1}{\sqrt{\beta}L}$. Hence, $T_N = \sum_{k=1}^N \tau_k \ge \frac{N\rho}{\sqrt{\beta}L}=N\rho\tau_0$. In particular, the bound~\eqref{AMP} in Theorem~\ref{thm:MP} reduces to, letting $\sigma_0=\beta\tau_0=\sqrt{\beta}/L$:
\begin{small}
\begin{equation}
\mathcal{G}_{x,y}(\hat x^N,\hat y^N) \le \frac{1}{\rho N}\left(\frac{1}{2\tau_0}\|x-x^1\|^2_2+ \frac{1}{\sigma_0}D_{\mathcal{Y}}(y,y^1)-\tau_1\theta_1(\mathcal{L}(\hat x^N,y)-\mathcal{L}(x^0,y))\right).
\end{equation}
\end{small}

{\bf The strongly convex case.} Consider now the case
$\gamma >0$. A first point is to show that $\tau_k\ge \frac{\rho}{\sqrt{\beta_k}L}$: this is guaranteed only
if the linesearch procedure is active, but since $\beta_k$
is decreased at each iteration one could fall into a situation
where this does not happen.

Let us assume that for some $k\ge 1$, $\tau_k\ge\frac{\rho}{\sqrt{\beta_k}L}$,
$\theta_k=\tau_k/\tau_{k-1}$, $\beta_k=\beta_{k-1}/(1+\gamma\tau_{k-1}\beta_{k-1})$,
and consider the next step. Then either the linesearch terminated after
more than one iteration and one must have $\tau_{k+1}\ge\frac{\rho}{\sqrt{\beta_{k+1}}L}$,
else one has
$\tau_{k+1} = \tau_k\sqrt{1+\theta_k}$. In that case,
\begin{multline*}
  \tau_{k+1} =\tau_k\sqrt{1+\theta_k} \ge\frac{\rho\sqrt{1+\theta_k}}{\sqrt{\beta_{k}}L}
  \ge\frac{\rho}{\sqrt{\beta_{k+1}}L}
  \Leftrightarrow
  \\
  (1+\theta_{k})\frac{\beta_k}{1+\gamma\tau_k\beta_k}\ge\beta_k
  \Leftrightarrow
  \theta_k=\frac{\tau_k}{\tau_{k-1}} \ge \gamma\tau_k\beta_k
  \Leftrightarrow
  1\ge \gamma\tau_{k-1}\beta_k = \frac{\gamma\tau_{k-1}\beta_{k-1}}{1+\gamma\tau_{k-1}\beta_{k-1}}
\end{multline*}
which is true. Hence, one will have that $\tau_k\ge \rho/(\sqrt{\beta_k} L)$
for all $k\ge 1$ provided $\tau_1\ge \rho/\sqrt{\beta_1}L$,
($\theta_1$, $\beta_1$ being defined from $\tau_1$, $\tau_0$, $\beta_0$
as in the algorithm).
A choice which guarantees that the induction holds is as follows:
\begin{equation}
  \beta_0>0,\quad \tau_0 = \frac{1}{\sqrt{\beta_0}L},\quad
  \theta_0 = \frac{\gamma\sqrt{\beta_0}}{L}.
\end{equation}
Then, one has
$\beta_1 = \beta_0/(1+\gamma\tau_0\beta_0)=\beta_0/(1+\gamma\sqrt{\beta_0}/L)
=\beta_0/(1+\theta_0)$ so that,
at the beginning of the linesearch for the
first step, $\tau_1 = \tau_0\sqrt{1+\theta_0}=\sqrt{(1+\theta_0)/\beta_0}/L
=1/(\sqrt{\beta_1}L)$ and in particular the first linesearch
will terminate immediately. The induction then guarantees
that for all $k\ge 2$, $\tau_k\ge\rho/(\sqrt{\beta_k}L)$.

Hence, one has for all $k\ge 1$:
\[
  \beta_k = \frac{\beta_{k-1}}{1+\gamma\beta_{k-1}\tau_{k-1}}
  \le \frac{\beta_{k-1}}{1+\frac{\gamma\rho}{L}\sqrt{\beta_{k-1}}}.
\]
One can write, for any $\alpha\in [0,1]$:
\[
  \frac{1}{\beta_k}\ge \frac{1}{\beta_{k-1}} + \frac{\gamma\rho}{L}\frac{1}{\sqrt{\beta_{k-1}}}
  = \left(\frac{1}{\sqrt{\beta_{k-1}}}+\frac{\alpha\gamma\rho}{2L}\right)^2
  +\frac{\gamma\rho}{L} \left(\frac{1-\alpha}{\sqrt{\beta_{k-1}}}
  -\frac{\alpha^2\gamma\rho}{4L}\right).
\]
Choosing to simplify $\alpha=.5$ one deduces that
$\beta_k\le 1/(1/\sqrt{\beta_l}+\gamma\rho(k-l)/(4L))^2$ provided
$\beta_l\le 64L^2/(\gamma\rho)^2$. Notice, though, that one also has,
by induction (assuming $L/(\gamma\rho)\ge 1$):
\[
  \beta_k \le \frac{L}{\gamma\rho}\sqrt{\beta_{k-1}}
  \le \frac{L^2}{\gamma^2\rho^2} \beta_0^{\frac{1}{2^k}}
\]
hence $\beta_l\le 64L^2/(\gamma\rho)^2$ as soon as
$l\ge \log_2\log_2\beta_0 -2$, which is very little whatever the choice
of $\beta_0$ ($l<8$ if $\beta_0\le 10^{300}$). So for ``reasonable'' choices
of $\beta_0$, and $k$ large enough, one has $\beta_k\le 16L^2/(\gamma\rho)^2 (1/k)^2$, $\tau_k\ge \gamma\rho^2(4L^2) k$, and $T_N\ge\gamma\rho^2 / (16L^2) N^2$.

\subsection{Termination of the linesearch} \label{termination}
In this section we prove that condition in step~\ref{MPrule} of Algorithm \ref{Alg:PDls} is always reached, and the number of iterations of the inner-loop does not affect the convergence rate. 

If $\tau_k \le \frac{1}{\sqrt{\beta_k}L}$ we have,
recalling that 
$L= \|K\|_{2\rightarrow\mathcal{Y}^*}=\|K^*\|_{\mathcal{Y}\rightarrow 2}
=\sup_{\|x\|_2\le 1,\|y\|_{\mathcal{Y}}\le 1} \langle x,K^*y\rangle$:
\begin{align*}
\tau_k\langle \bar x^k-x^{k+1}, K^*y^{k+1}-K^*y^k\rangle & \le \tau_k\|x^{k+1}-\bar x^k\|_{2} \|K^*y^{k+1}-K^*y^k\|_{2} \\
& \le \frac{1}{2}\|x^{k+1}-\bar x^k\|^2_{2}+ \frac{\tau_k^2}{2}L^2\|y^{k+1}-y^k\|_\mathcal{Y}^2 \\
& \le \frac{1}{2}\|x^{k+1}-\bar x^k\|^2_{2}+ \frac{1}{2\beta_k}\|y^{k+1}-y^k\|_\mathcal{Y}^2 \\
& \le \frac{1}{2}\|x^{k+1}-\bar x^k\|^2_2+ \frac{1}{\beta_k}D(y^{k+1},y^k) 
\end{align*}
Therefore, the condition is reached as soon as $\tau_k \le \frac{1}{\sqrt{\beta_k}L}$ which eventually will occur since $\tau_k$ decrease by a factor $\rho<1$ at each iteration of the inner-loop. We estimate, following~\cite{malitsky2018first}, the extra computing time due to
  the inner loops, up to the $k$th outer loop.

By construction of the algorithm $\theta_k \le \sqrt{1+\theta_{k-1}}$, hence by induction we have that $\theta_k\le \bar \theta := \max\left\{ \frac{1+\sqrt{5}}{2},\theta_0\right\}$. Fix iteration $k$. After $i_k$ iterations of the inner-loop we have $\tau_k(i_k) = \tau_{k-1}\sqrt{1+\theta_{k-1}}\rho^{i_k}$ and $\theta_k(i_k) = \sqrt{1+\theta_{k-1}}\rho^{i_k}\le \bar \theta \rho^{i_k}$. Let $\tau_0$ such that the condition in step~\ref{MPrule} is true (for instance $\tau_0 = \frac{1}{\sqrt{\beta_0}L}$), then:
$$\bar \theta \rho^{i_k}\ge \theta_k \ge \frac{\tau_0\rho}{\tau_{k-1}} = \frac{\rho}{\prod_{j=1}^{k-1}\theta_j} \ge \frac{\rho^{k+1}}{\bar \theta^k \rho^{\sum_{j=1}^{k-1}i_j}},$$
so that
\begin{equation}
  \sum_{j=1}^{k} i_j \le (k+1)\left(1+\frac{\ln \bar \theta}{|\ln \rho|}\right).
\end{equation}

\subsection{Proof of Lemma~\ref{thm:halfC}}\label{sec:halfC}

We assume as in the statement that $\min_{i} C_{i,j} = 0$ for all $j = 1,...,n$, and $\min_{j} C_{i,j} = 0$ for all $i=0,...,n$. (Otherwise we remove the minimum on each row and then  on each column of $C$ without
changing the problem, up to a constant.)

Let $(u,v)$ a dual solution of the saddle point problem. Define $\tilde{u} = u-\onen\min_i u_i$ and $\tilde{v} = v+\onen\min_i u_i$. The new dual solution has the same value, $\tilde{u}\ge 0$ and there exist $i_0$ such that $\tilde{u}_{i_0} = 0$. The dual constraints are $\tilde{u}_i+\tilde{v}_j\le C_{ij}$ for all $i,j$. By optimality of $(\tilde{u},\tilde{v})$ we have $\tilde{u}_i = \min_{j}(C_{i,j}-\tilde{v}_j)$ for each $i$, and in particular there exist $j_0$ such that $\tilde{u}_{i_0}+\tilde{v}_{j_0}= \tilde{v}_{j_0}= C_{i_0,j_0}$.
Symmetrically, we have $\tilde{v}_j = \min_i(C_{i,j}-\tilde{u}_i)\le \min_i C_{i,j} = 0$ since $\tilde{u}_i\ge 0$.
In particular $\tilde{v}_{j_0}=C_{i_0,j_0}=0$.

As a consequence, for each $i$, $\tilde u_i\le C_{i,j_0}-\tilde v_{j_0}\le \|C\|$,
and for each $j$, $\exists i$ with $\tilde{v}_j = C_{i,j}-\tilde{u}_{\color{red} i}
\ge C_{i,j}-\|C\|\ge -\|C\|$.
In brief, we have $0\le \tilde{u}_i \le \|C\|$ and $-\|C\| \le \tilde{v}_j \le 0$. Defining $u^* = \tilde{u}-\onen\frac{\|C\|}{2}$ and  $v^* = \tilde{v}+\onen\frac{\|C\|}{2}$ we have the first result. 

For the second part, take $\mu = (1,0,...,0)^\top$ and $\nu = (0,...,0,1)^\top$, and $C = \mu \otimes \nu$. The unique optimal transport in this case is $X_{1,n} = 1$ and $X_{ij} = 0$ otherwise. Therefore, for any dual solution $(u^*,v^*)$ we have $C_{1,n} = 1 = u_1^*+v_n^*$, hence $\|(u^*,v^*)\| \ge 1/2 = \frac{\|C\|}{2}$. 

\subsection{Proof of Lemma~\ref{thm:approx2}}\label{sec:approx2}

We show that Algorithm~2, a.k.a.~``\textsc{Round}'', in~\cite{altschuler2017near}
satisfies a better error bound than proved in the original paper,
at least if its input is already in the matrix unit simplex.
We consider $X\in {\Delta}$: hence one has $\sum_{i,j}X_{i,j}=1$,
and $\sum_i(X\onen -\mu)_i=0$, $\sum_j(X^\top\onen-\nu)_j=0$.

This algorithm outputs $Y=\mathcal{A}X$ as follows: First, one lets:
\[
  X'_{i,j} = \min\{1,\tfrac{\mu_i}{(X\onen)_i}\} X_{i,j}
  = \begin{cases}
    X_{i,j} & \text{ if } \mu_i \ge (X\onen)_i\\
    \tfrac{\mu_i}{(X\onen)_i} X_{i,j} & \text{ else}
  \end{cases}
\]
and in particular $X'\le X$ and $X'\onen\le \mu$. Then:
\[
  X''_{i,j} = X'_{i,j} \min\{1,\tfrac{\nu_j}{((X')^\top\onen)_j}\}
  = \begin{cases}
    X'_{i,j} & \text{ if } \nu_j \ge ((X')^\top\onen)_j\\
    \tfrac{\nu_j}{((X')^\top\onen)_j} X'_{i,j} & \text{ else}
  \end{cases}
\]
which is such that $X''\le X'$ and $X''^\top\onen\le\nu$.
In particular since $\mu\ge X'\onen\ge X''\onen$, $\|\mu-X''\onen\|_1=\sum_i(\mu_i-\sum_j X''_{i,j})= 1-\sum_{i,j}X''_{i,j}
= \sum_j \left(\nu_j-\sum_iX''_{i,j}\right)$.
Eventually the output is given by:
\[
  Y_{i,j} = X''_{i,j}+\frac{\mu_i-(X''\onen)_i}{\|\mu-X''\onen\|_1}(\nu_j-((X'')^\top\onen)_j).
\]
and we see that $Y\onen=\mu$, $Y^\top\onen=\nu$.
One has, then:
\[
  Y-X = Y-X'' - (X'-X'') - (X-X').
\]
Let $J^+=\{j:\nu_j\ge ((X')^\top\onen)_j$ and $J^-=\{1,\dots,N\}\setminus J^+$.
Then if $j\in J^+$, $X''_{i,j}=X'_{i,j}$ and
\begin{align*}
  Y_{i,j}-X_{i,j}
  &=\frac{\mu_i-(X''\onen)_i}{\|\mu-X''\onen\|_1}(\nu_j-((X'')^\top\onen)_j) -(X_{i,j}-X'_{i,j})
  \\
  &=\frac{\mu_i-(X''\onen)_i}{\|\mu-X''\onen\|_1}(\nu_j-((X')^\top\onen)_j) -(X_{i,j}-X'_{i,j})
  \\
  &=\frac{\mu_i-(X''\onen)_i}{\|\mu-X''\onen\|_1}(\nu_j-(X^\top\onen)_j)
       +\frac{\mu_i-(X''\onen)_i}{\|\mu-X''\onen\|_1}((X-X')^\top\onen)_j
        -(X_{i,j}-X'_{i,j})
\end{align*}
and
\[
  \sum_i |Y_{i,j}-X_{i,j}| \le |\nu_j-(X^\top\onen)_j| + 2((X-X')^\top\onen)_j.
\]
If $j\in J^-$, then $((X'')^\top\onen)_j=\nu_j$ hence $Y_{i,j}=X''_{i,j}$, so that:
\begin{align*}
  Y_{i,j}-X_{i,j}
  &= -(X'_{i,j}-X''_{i,j})-(X_{i,j}-X'_{i,j})
  \\
  & = -\frac{((X')^\top\onen)_j-\nu_j}{((X')^\top\onen)_j} X'_{i,j}
    -(X_{i,j}-X'_{i,j}) \\ &
                             = - \frac{X'_{i,j}}{((X')^\top\onen)_j}((X^\top\onen)_j-\nu_j)
       - \frac{X'_{i,j}}{((X')^\top\onen)_j}((X'-X)^\top\onen)_j
    - (X_{i,j}-X'_{i,j}).
\end{align*}
and again
\[
  \sum_i |Y_{i,j}-X_{i,j}| \le (X^\top\onen)_j-\nu_j + 2((X-X')^\top\onen)_j.
\]
Hence,
\[
  \|Y-X\|_1 \le \|X^\top\onen -\nu\|_1 + 2\|X-X'\|_1.
\]
Eventually, one has:
\[
  \|X-X'\|_1 = \frac{1}{2}\|X\onen-\mu\|_1
\]
since for all $i,j$,
\[
  X_{i,j}-X'_{i,j} = \frac{((X\onen)_i-\mu_i)^+}{(X\onen)_i} X_{i,j}
\]
and $\sum_{i} (X\onen)_i-\mu_i=0\Rightarrow
\sum_{i}((X\onen)_i-\mu_i)^+=\sum_{i}((X\onen)_i-\mu_i)^-
=\frac{1}{2}\sum_{i}|(X\onen)_i-\mu_i)|$.
It follows that
\begin{equation}
  \|Y-X\|_1 \le \|X^\top\onen -\nu\|_1 + \|X\onen-\mu\|_1.
\end{equation}

\end{document}